# OPTIMAL SCALING FOR PARTIALLY UPDATING MCMC ALGORITHMS

By Peter Neal and Gareth Roberts

*University of Manchester and Lancaster University*

In this paper we shall consider optimal scaling problems for high-dimensional Metropolis–Hastings algorithms where updates can be chosen to be lower dimensional than the target density itself. We find that the optimal scaling rule for the Metropolis algorithm, which tunes the overall algorithm acceptance rate to be 0.234, holds for the so-called Metropolis-within-Gibbs algorithm as well. Furthermore, the optimal efficiency obtainable is independent of the dimensionality of the update rule. This has important implications for the MCMC practitioner since high-dimensional updates are generally computationally more demanding, so that lower-dimensional updates are therefore to be preferred. Similar results with rather different conclusions are given for so-called Langevin updates. In this case, it is found that high-dimensional updates are frequently most efficient, even taking into account computing costs.

**1. Introduction.** There exist large classes of Markov chain Monte Carlo (MCMC) algorithms for exploring high-dimensional (target) distributions. All methods construct Markov chains with invariant distribution given by the target distribution of interest. However, for the purposes of maximizing the efficiency of the algorithm for Monte Carlo use, it is imperative to design algorithms which give rise to Markov chains which mix sufficiently rapidly. Since all Metropolis–Hastings algorithms require the specification of a proposal distribution, these implementational questions can all be phrased in terms of proposal choice. This paper is about two of these choices: the scaling and dimensionality of the proposal. We shall work throughout with continuous distributions, although it is envisaged that more general distributions might be amenable to similar study.

One important decision the MCMC user has to make in a $d$-dimensional problem concerns the dimensionality of the proposed jump. For instance, two extreme types of algorithm are the following: propose a fully $d$-dimensional update of the current state (according to a density with a density with respect to $d$-dimensional Lebesgue measure) and accept or reject according to the Metropolis–Hastings acceptance probabilities; or, for each of the $d$ components in turn, update that component conditional on all the others according to some Markov chain which preserves the appropriate conditional distribution. The most widely used example is the $d$-dimensional Metropolis algorithm, in one extreme, and the Gibbs sampler or some kind of "Metropolis-within-Gibbs" scheme in the other. In between these two options, there lie many intermediate strategies. An important question is whether any general statements can be made about algorithm choice in this context, leading to practical advice for MCMC practitioners.

In this paper we concentrate on two types of algorithm: Metropolis and Metropolis adjusted Langevin algorithms (MALA). We consider strategies which update a fixed proportion, $c$, of components at each iteration, and consider the efficiency of the algorithms constructed asymptotically as $d \to \infty$. In order to do this, we shall extend the methodology developed in [6, 7] to our context. The analysis produces clear cut results which suggest that, while full-dimensional Langevin updates are worthwhile, full-dimensional Metropolis ones are asymptotically no better than smaller dimensional updating schemes, so that the possible extra computational overhead associated with their implementation always leads to their being suboptimal in practice. All this is initially done in the context of target densities consisting of independent components, and this leads naturally to the question of whether this simple picture is altered in any way in the presence of dependence. Although this is difficult to explore in full generality, we do later consider this problem in the context of a class of Gaussian dependent target distributions where explicit results can be shown, and where the conclusions from the independent component case remain valid.

It is now well recognized that highly correlated target distributions lead to slow mixing for updating schemes where $c < 1$ (see, e.g., [5, 9]). However, it is also known that spherically symmetric proposal distributions in $d$-dimensions on highly correlated target densities can lead to slow mixing since the proposal distribution is inappropriately shaped to explore the target (see [8]). So for highly correlated target distributions, both high and small dimensional updating strategies perform poorly. We shall explore these two competing algorithms in a Gaussian context where explicit calculations are possible. Our work shows that, for $c > 0$, for the Metropolis algorithm, these two slowing down effects are the same. In particular, this implies that the commonly used strategy of getting round high correlation problems by block updating using Metropolis has no justification. In contrast, for MALA full dimensional updating, $c = 1$, is shown to be optimal.

The paper is structured as follows. In Section 2 we outline the MCMC setup. In Sections 3 and 4 we tackle the problem of scaling the variance of the proposal distribution for RWM-within-Gibbs (Random walk Metropolis-within-Gibbs) and MALA-within-Gibbs (Metropolis adjusted Langevin-within-Gibbs), respectively. The approach taken is similar to that used for the full RWM/MALA algorithms, by obtaining weak convergence to an appropriate Langevin diffusion as the dimension of the state space, $d$ converges to infinity. The results of Sections 3 and 4 are proved for a sequence of $d$-dimensional product densities of the form

$$\pi_d(\mathbf{x}^d) = \prod_{i=1}^{d} f(x_i^d) \tag{1.1}$$

for some suitably smooth probability density $f(\cdot)$. In both Sections 3 and 4, for each fixed, one-dimensional component of $\{\mathbf{X}^d; d \geq 1\}$, the one-dimensional process converges weakly to an appropriate Langevin diffusion. The aim therefore is to scale the proposal variances so as to maximize the speed of the limiting Langevin diffusion. Since each of the components of $\{\mathbf{X}^d; d \geq 1\}$ are independent and identically distributed, we shall prove the results for $\{X_1^d; d \geq 1\}$.

However, it is at least plausible that the picture will be very different when considering dependent densities. However, theoretical analysis in the limiting case where results can be obtained and in simulations for more general cases, we find that the general conclusions which can be derived for densities of the form (1.1) extend some way toward dependent densities. To this end, in Section 5, we consider RWM/MALA-within-Gibbs for the exchangeable normal $\mathbf{X}^d \sim N(\mathbf{0}, \Sigma_\rho^d)$, where $\sigma_{ii}^d = 1$, $1 \leq i \leq d$, and $\sigma_{ij}^d = \sigma_{ji}^d = \rho$, $1 \leq i < j \leq d$. (Throughout the paper, we adopt the notation that $\Sigma$ will be used for variance matrices, while elements of matrices will be denoted by $\sigma$, both conventions using appropriate sub- and super-scripts.)

All the proofs of the theorems in Sections 3–5 are given in the Appendix. Then in Section 6 with the aid of a simulation study we demonstrate that the asymptotic results are practically useful for finite $d$, namely, $d \geq 10$.

## 2. Algorithms and preliminaries.

For RWM/MALA, we are interested in $(d, \sigma_d^2)$, the dimension of the state space, $d$, and the proposal variance $\sigma_d^2$, where the proposal for the $i$th component is given by

$$Y_i^d = x_i^d + \sigma_d Z_i, \qquad 1 \leq i \leq d, \qquad \text{RWM},$$

$$Y_i^d = x_i^d + \sigma_d Z_i + \frac{\sigma_d^2}{2} \frac{\partial}{\partial x_i} \log \pi_d(\mathbf{x}^d), \qquad 1 \leq i \leq d, \qquad \text{MALA}$$

and the $\{Z_i\}$'s are independent and identically distributed according to $Z \sim N(0,1)$. For both RWM and MALA, the maximum speed of the diffusion can be obtained by taking the proposal variance to be of the form $\sigma_d^2 = l^2 d^{-s}$ for some $l > 0$ and $s > 0$. (For RWM, $s = 1$ and for MALA, $s = \frac{1}{3}$.)

Now for RWM/MALA-within-Gibbs, the basic idea is to choose $dc_d$ components at random at each iteration, attempting to update them jointly according to the RWM/MALA mechanism, respectively. We sometimes write $\sigma_d^2 = \sigma_{d,c_d}^2$, where $c_d$ represents the proportion of components updated at each iteration. Thus, the two algorithms propose new values as follows:

$$
\begin{aligned}
&Y_i^d = x_i^d + \chi_i^d \sigma_{d,c_d} Z_i, \qquad 1 \le i \le d, \qquad \text{RWM-within-Gibbs},\\
&Y_i^d = x_i^d + \chi_i^d \left\{ \sigma_{d,c_d} Z_i + \frac{\sigma_{d,c_d}^2}{2} \frac{\partial}{\partial x_i} \log \pi_d(\mathbf{x}^d) \right\},\\
&\qquad\qquad 1 \le i \le d, \qquad \text{MALA-within-Gibbs},
\end{aligned}
\tag{2.1}
$$

where the $\{Z_i\}$'s are independent and identically distributed according to $Z \sim N(0,1)$ and the $\{\chi_i^d\}$ are chosen as follows. Independently of the $Z_i$'s, we select at random a subset $A$, say, of size $dc_d$ from $\{1, 2, \ldots, d\}$, setting $\chi_i^d = 1$ if $i \in A$, and $\chi_i^d = 0$ otherwise. The proposal $\mathbf{Y}^d$ is then accepted according to the usual Metropolis–Hastings acceptance probability:

$$
\alpha_d^{c_d}(\mathbf{x}^d, \mathbf{Y}^d) = 1 \wedge \frac{\pi_d(\mathbf{Y}^d) q(\mathbf{Y}^d, \mathbf{x}^d)}{\pi_d(\mathbf{x}^d) q(\mathbf{x}^d, \mathbf{Y}^d)},
\tag{2.2}
$$

where $q(\cdot,\cdot)$ is the proposal density. Otherwise, we set $\mathbf{X}_m^d = \mathbf{X}_{m-1}^d$.

In both cases, the algorithms simulate Markov chains which are reversible with respect to $\pi_d$, and can be easily shown to be $\pi_d$-irreducible and aperiodic. Therefore, both algorithms will converge in total variation distance to $\pi_d$. However, here we shall investigate optimization of the algorithms for rapid convergence. To find a manageable framework for assessing optimality, Roberts, Gelman and Gilks [6] introduce the notion of the average acceptance rate which measures the steady state proportion of accepted proposals for the algorithm, and which can be shown to be closely connected with the notion of algorithm efficiency and optimality. Specifically, we define

$$
a_d^{c_d}(l) = \mathbb{E}_{\pi_d}[\alpha_d^{c_d}(\mathbf{X}^d, \mathbf{Y}^d)] = \mathbb{E}_{\pi_d}\left[1 \wedge \frac{\pi_d(\mathbf{Y}^d) q(\mathbf{Y}^d, \mathbf{X}^d)}{\pi_d(\mathbf{X}^d) q(\mathbf{X}^d, \mathbf{Y}^d)}\right],
\tag{2.3}
$$

where $\sigma_{d,c_d}^2 = l^2 d^{-s}$, $\mathbf{X}^d \sim \pi_d$ and $\mathbf{Y}^d$ represents the subsequent proposal random variable. Thus, $a_d^{c_d}(l)$ is the $\pi_d$-average acceptance rate of the above algorithms where we update a proportion $c_d$ of the $d$ components in each iteration. We adopt the general notational convention that, for any $d$-dimensional stochastic process $\mathbf{W}_{\cdot,\cdot}^d$, we shall write $W_{t,i}^d$ for the value of its $i$th component at time $t$.

Our aim in this paper is to consider the optimization [in $(c_d, \sigma^2_{d,c_d})$] of the algorithms speed of convergence. For convenience (although to some extent this assumption can be relaxed), we shall assume that $c_d \to c$ as $d \to \infty$ for some $0 < c \leq 1$. It turns out to be both convenient and practical to express many of the optimality solutions in terms of acceptance rate criteria.

**3. RWM-within-Gibbs for IID product densities.** We shall first consider the RWM algorithm applied initially to a simple IID form target density. This allows us to obtain explicit asymptotic results for optimal high-dimensional algorithms. The results of this section can be seen as an extension of the results of Theorems 1.1 and 1.2 of [6] which considers the full-dimensional update case.

Let

$$\pi_d(\mathbf{x}^d) = \prod_{i=1}^{d} f(x_i^d) = \prod_{i=1}^{d} \exp\{g(x_i^d)\} \tag{3.1}$$

be a $d$-dimensional product density with respect to Lebesgue measure. Let the proposal standard deviation $\sigma_d = \frac{l}{\sqrt{d-1}}$ for some $l > 0$.

For $d \geq 1$, let $\mathbf{U}_t^d = (X^d_{[dt],1}, X^d_{[dt],2}, \ldots, X^d_{[dt],d})$, and so, $U^d_{t,i} = X^d_{[dt],i}$, $1 \leq i \leq d$. Let $U_t^d = U^d_{t,1}$.

THEOREM 3.1. *Suppose that $f$ is positive, $C^3$ (a three-times differentiable function with continuous third derivative) and that $(\log f)' = g'$ is Lipschitz. Suppose also that, $c_d \to c$, as $d \to \infty$, for some $0 < c \leq 1$,*

$$\mathbb{E}_f\left[\left(\frac{f'(X)}{f(X)}\right)^8\right] < \infty \tag{3.2}$$

*and*

$$\mathbb{E}_f\left[\left(\frac{f''(X)}{f(X)}\right)^4\right] < \infty. \tag{3.3}$$

*Let $\mathbf{X}_0^\infty = (X^1_{0,1}, X^2_{0,2}, \ldots)$ be such that all of its components are distributed according to $f$ and assume that $X^j_{0,i} = X^i_{0,i}$ for all $i \leq j$. Then, as $d \to \infty$,*

$$U^d \Rightarrow U, \tag{3.4}$$

*where $U_0$ is distributed according to $f$ and $U$ satisfies the Langevin SDE*

$$dU_t = (h_c(l))^{1/2}\, dB_t + \tfrac{1}{2} h_c(l) g'(U_t)\, dt \tag{3.5}$$

*and*

$$h_c(l) = 2cl^2 \Phi\left(-\frac{l\sqrt{cI}}{2}\right),$$

*with* $\Phi$ *being the standard normal cumulative c.d.f and*

$$I \equiv \mathbb{E}_f\left[\left(\frac{f'(X)}{f(X)}\right)^2\right] \equiv \mathbb{E}_g[g'(X)^2].$$

The following corollary holds.

COROLLARY 3.2. *Let* $c_d \to c$, *as* $d \to \infty$, *for some* $0 < c \leq 1$. *Then:*

(i) $\lim_{d\to\infty} a_d^{c_d}(l) = a^c(l) \stackrel{def}{=} 2\Phi(-\frac{l\sqrt{cI}}{2})$.

(ii) *Let* $\hat{l}$ *be the unique value of* $l$ *which maximizes* $h_1(l) = 2l^2\Phi(-\frac{l\sqrt{I}}{2})$ *on* $[0,\infty)$, *and let* $\hat{l}_c$ *be the unique value of* $l$ *which maximizes* $h_c(l)$ *on* $[0,\infty)$. *Then* $\hat{l}_c = c^{-1/2}\hat{l}$ *and* $h_c(\hat{l}_c) = h_1(\hat{l})$.

(iii) *For all* $0 < c \leq 1$, *the optimal acceptance rate* $a^c(\hat{l}_c) = 0.234$ *(to three decimal places).*

Though these results involve fairly technical mathematical statements, they yield a very simple practical conclusion. Optimal efficiency obtainable for a given $c$ does not depend on $c$ at all. Now, in practice, computational overheads associated with one iteration of the algorithm are nondecreasing as a function of $c$, so that, in practice, smaller values of $c$ should be preferred. Therefore, for RWM, using high-dimensional update steps does not make any sense.

It is, of course, important to see how these conclusions extend to more general target densities and, in particular, ones which exhibit dependence structure. Some theory and related simulation studies in Sections 5 and 6, respectively, will demonstrate that these findings extend considerably beyond the rigorous but restrictive set up of Theorem 3.1.

**4. MALA-within-Gibbs for IID product densities.** We now turn our attentions to MALA-within-Gibbs. We again consider a sequence of probability densities $\pi_d$ of the form given in (3.1). We follow [7] in making the following assumptions. We assume that $\mathbf{X}_0^d$ is distributed according to the stationary measure $\pi_d$, $g$ is an eight times continuously differentiable function with derivatives $g^{(i)}$ satisfying

$$|g(x)|, |g^{(i)}(x)| \leq C(1 + |x|^K), \tag{4.1}$$

$1 \leq i \leq 8$, for some $C, K > 0$, and that

$$\int_{\mathbb{R}} x^k f(x)\, dx < \infty, \qquad k = 1, 2, \ldots. \tag{4.2}$$

Finally, we assume that $g'$ is Lipschitz. This ensures that $\{\mathbf{X}_t\}$ is nonexplosive (see, e.g., [12], Chapter V, Theorem 52.1).

Let $\{J_t\}$ be a Poisson process with rate $d^{1/3}$ and let $\Gamma^d = \{\Gamma_t^d\}_{t\geq 0}$ be the $d$-dimensional jump process defined by $\Gamma_t^d = \mathbf{X}_{J_t}^d$, where we take $\sigma_d^2 = l^2 d^{-1/3}$ with $l$ an arbitrary constant.

We then have the following two theorems which are extensions of [7], Theorems 1 and 2.

THEOREM 4.1. *Suppose that $c_d \to c$, as $d \to \infty$, for some $0 < c \leq 1$. We have that*

$$\lim_{d\to\infty} \{a_d^{c_d}(l)\} = a^c(l) = 2\Phi\left(-\frac{\sqrt{c}Kl^3}{2}\right),$$

*with $K^2 = \mathbb{E}[\frac{5g'''(X)^2 - 3g''(X)^3}{48}] > 0$.*

THEOREM 4.2. *Suppose that $c_d \to c$ as $d \to \infty$ for some $0 < c \leq 1$. Let $\{U^d\}_{t\geq 0}$ be the process corresponding to the first component of $\Gamma^d$. Then, as $d \to \infty$, the process $U^d$ converges weakly (in the Skorokhod topology) to the Langevin diffusion $U$ defined by*

$$dU_t = h_c(l)^{1/2}\, dB_t + \tfrac{1}{2} h_c(l) g'(U_t)\, dt,$$

*where $h_c(l) = 2cl^2\Phi(-\frac{\sqrt{c}l^3K}{2})$ is the speed of the limiting diffusion.*

The most important consequence of Theorems 4.1 and 4.2 is the following corollary.

COROLLARY 4.3. *Let $c_d \to c$, as $d \to \infty$, for some $0 < c \leq 1$. Then:*

(i) *Let $\hat{l}$ be the unique value of $l$ which maximizes $h_1(l) = 2l^2\Phi(-\frac{l^3K}{2})$ on $[0,\infty)$, and let $\hat{l}_c$ be the unique value of $l$ which maximizes $h_c(l)$ on $[0,\infty)$. Then $\hat{l}_c = c^{-1/6}\hat{l}$ and $h_c(\hat{l}_c) = c^{2/3} h_1(\hat{l})$.*

(ii) *For all $0 < c \leq 1$, the optimal acceptance rate $a^c(\hat{l}_c) = 0.574$ (to three decimal places).*

Thus, in stark contrast to the RWM case, it is optimal to update all components at once for MALA. The story is somewhat more complicated in the case where computational overheads are taken into account. For instance, it is common for the computational costs of implementing MALA-within-Gibbs to be approximately $d(a + bc)$ for constants $a$ and $b$. To see this, note that the algorithm's computational cost is often dominated by two operations: the calculation of the various derivatives needed to propose a new value, and the evaluation of $\pi$ at the proposed new value. The first of these operations involves a $cd$-dimensional update and typically takes a time which is order $cd$, while the second involves evaluating a $d$-dimensional function which we would expect to be at least of order $d$. (Although, in some important special cases, target density ratios might be computed more efficiently than this.) In this case the overall efficiency is obtained by maximizing

$$\frac{c^{2/3}}{a+bc}.$$

This expression is maximized at $1 \wedge 2a/b$. Therefore, it is conceivable for full dimensional updates to be optimal even when computational costs are taken into account. In any case, the optimal proportion will be some value $x^* \in (0,1]$.

**5. RWM/MALA-within-Gibbs on dependent target distributions.** We are now interested in the extent to which the results of the last two sections can be extended to the case where the $d$ components are dependent. It is difficult to get general results, but certain important special cases can be examined explicitly, yielding interesting results which imply (essentially) that the extent by which the dependence structure affects the mixing properties of the chain (RWM-within-Gibbs or MALA-within-Gibbs) is independent of $c$. The most tractable special case is the Gaussian target distribution. However, in Section 6, we shall also include some simulations in other cases to show that the above statement holds well beyond the cases for which rigorous mathematical results can be proved.

We begin with RWM-within-Gibbs and consider the optimal scaling problem of the variance of the proposal distribution for a target distribution consisting of exchangeable normal components. Specifically, $\mathbf{X}^d \sim N_d(\mathbf{0}, \Sigma_\rho^d)$, where $\sigma_{ii}^d = 1$, $1 \le i \le d$, and $\sigma_{ij}^d = \rho$, $i \ne j$, for some $0 < \rho < 1$. Therefore, we have that

$$\begin{aligned} \pi_d(\mathbf{x}^d) &= (2\pi)^d \det|\Sigma_\rho^d|^{-1/2} \\ &\quad \times \exp\left(-\frac{1}{2}\left(\frac{1}{1-\rho}\sum_{i=1}^d (x_i^d)^2 + \theta_d \sum_{i=1}^d \sum_{j=1}^d x_i^d x_j^d\right)\right) \qquad (5.1) \\ &= (2\pi)^d \det|\Sigma_\rho^d|^{-1/2} \exp\left(-\frac{1}{2} j_d(\mathbf{x}^d)\right), \qquad \text{say,} \end{aligned}$$

where

$$\theta_d = \frac{-\rho}{1+(d-2)\rho-(d-1)\rho^2}$$

and

$$j_d(\mathbf{x}^d) = \frac{1}{1-\rho}\sum_{i=1}^d (x_i^d)^2 + \theta_d \sum_{i=1}^d \sum_{j=1}^d x_i^d x_j^d.$$

For $d \geq 1$, $\hat{\mathbf{U}}_t^d = (X^d_{[dt],1}, X^d_{[dt],2}, \ldots, X^d_{[dt],d})$. Let $\mathbf{U}_t^d = (U^d_{t,1}, U^d_{t,2}, U^d_{t,3})$ be such that $U^d_{t,1} = \hat{U}^d_{t,1}$, $U^d_{t,2} = \hat{U}^d_{t,2}$ and $U^d_{t,3} = \frac{1}{d-2}\sum_{i=3}^d \hat{U}_{t,i}$.

Now the proposal $\mathbf{Y}^d$ is given by

$$Y_i^d = x_i^d + \sigma_d \chi_i^d Z_i, \qquad 1 \leq i \leq d,$$

where the $Z_i$ and $\chi_i^d$ $(1 \leq i \leq d)$ are defined as before and $\sigma_d = \frac{l}{\sqrt{d-2}}$ for some constant $l$. [We use $(d-2)$ rather than $d$ or $(d-1)$ for simplicity in presentation of the results.]

In the dependent case, more care needs to be taken in constructing the sequence $\{\mathbf{X}_0^d; d \geq 1\}$. Let $\mathbf{X}_0^1 \sim N(0,1)$ [i.e., $\mathbf{X}_0^1$ is distributed according to $\pi_1(\cdot)$]. For $d \geq 2$ and $1 \leq i \leq d-1$, set $X^d_{0,i} = X^i_{0,i}$. Then iteratively define

$$X^d_{0,d} \sim N\left(\rho \frac{1}{d-1}\sum_{i=1}^{d-1} X^d_{0,i}, \frac{1}{d-1}(1 + (d-2)\rho - (d-1)\rho^2)\right).$$

Therefore, $\mathbf{X}_0^d$ is distributed according to $\pi_d(\cdot)$ and we can continue this process indefinitely to obtain $\mathbf{X}_0^\infty = (X^1_{0,1}, X^2_{0,2}, \ldots)$.

THEOREM 5.1. *Suppose that* $0 < \rho < 1$ *and that* $c_d \to c$, *as* $d \to \infty$, *for some* $0 < c \leq 1$. *Let* $\mathbf{X}_0^\infty = (X^1_{0,1}, X^2_{0,2}, \ldots)$ *be constructed as above. Let*

$$D_1 = \begin{pmatrix} 1 & \rho & \rho \\ \rho & 1 & \rho \\ \rho & \rho & \rho \end{pmatrix},$$

$$D_2 = \begin{pmatrix} \dfrac{1}{1-\rho} & 0 & -\dfrac{1}{1-\rho} \\ 0 & \dfrac{1}{1-\rho} & -\dfrac{1}{1-\rho} \\ -\dfrac{1}{1-\rho} & -\dfrac{1}{1-\rho} & \dfrac{1+\rho}{\rho(1-\rho)} \end{pmatrix},$$

$$D_3 = \begin{pmatrix} 1 & 0 & 0 \\ 0 & 1 & 0 \\ 0 & 0 & 0 \end{pmatrix}.$$

*Let* $\tilde{f}(\mathbf{u})$ *denote the probability density function of* $N(0, D_1)$. *Then, as* $d \to \infty$,

$$\mathbf{U}^d \Rightarrow \mathbf{U},$$

*where* $\mathbf{U}_0$ *is distributed according to* $\tilde{f}$ *and* $\mathbf{U}$ *satisfies the Langevin SDE*

$$d\mathbf{U}_t = (h_{c,\rho}(l))^{1/2} D_3\, d\mathbf{B}_t + h_{c,\rho}(l) D_3 \{\tfrac{1}{2}\nabla(-\tfrac{1}{2}\mathbf{U}_t^T D_2 \mathbf{U}_t)\}\, dt,$$

*where*

$$h_{c,\rho}(l) = 2cl^2 \Phi\left(-\frac{l}{2}\sqrt{\frac{c}{1-\rho}}\right).$$

Note that if we define $\tilde{I}_d = \mathbb{E}[(\frac{\partial}{\partial x_1} j_d(\mathbf{X}^d))^2]$ and $\tilde{I} = \frac{1}{1-\rho}$. Then $\tilde{I}_d \to \tilde{I}$ as $d \to \infty$ and $h_{c,\rho}(l) = 2cl^2\Phi(-\frac{l}{2}\sqrt{c\tilde{I}})$. Therefore, the speed of the limiting diffusion for exchangeable normal has the same form as that obtained for the IID product densities considered in Section 3.

As in (2.3), let $a_d^{c_d,\rho}(l)$ be the $\pi_d$-average acceptance rate of the above algorithm where $\mathbf{X}^d \sim N(\mathbf{0}, \Sigma_\rho^d)$, $\sigma_d = \frac{l}{\sqrt{d-2}}$ and we update a proportion $c_d$ of the $d$ components in each iteration. Then we have the following corollary.

Corollary 5.2. *Let $c_d \to c$, as $d \to \infty$, for some $0 < c \leq 1$. Then, for $0 < \rho < 1$:*

(i) $\lim_{d\to\infty} a_d^{c_d,\rho}(l) = a^{c,\rho}(l) \stackrel{def}{=} 2\Phi(-\frac{l}{2}\sqrt{\frac{c}{1-\rho}})$.

(ii) *Let $\hat{l}$ be the unique value of $l$ which maximizes $h_{1,0}(l) = 2l^2\Phi(-\frac{l}{2})$ on $[0,\infty)$, and let $\hat{l}_{c,\rho}$ be the unique value of $l$ which maximizes $h_{c,\rho}(l)$ on $[0,\infty)$. Then $\hat{l}_{c,\rho} = \sqrt{\frac{1-\rho}{c}}l$ and $h_{c,\rho}(\hat{l}_{c,\rho}) = (1-\rho)h_{1,0}(\hat{l})$.*

(iii) *For all $0 < c \leq 1$ and $0 < \rho < 1$, the optimal acceptance rate $a^{c,\rho}(\hat{l}_{c,\rho}) = 0.234$ (to three decimal places).*

Note that Corollary 5.2(ii) states that the cost incurred by having $\sigma_{ij}^d = \rho$, $i \neq j$, rather than $\sigma_{ij}^d = 0$, $i \neq j$, is to slow down the speed of the limiting diffusion by a factor of $1-\rho$, for all $0 < c \leq 1$. In other words, the cost incurred by the dependence between the components of $\mathbf{X}^d$ is independent of $c$. Furthermore, the optimal acceptance rate $a^{c,\rho}(\hat{l}_{c,\rho})$ is unaffected by the introduction of dependence. We shall study this further in the simulation study conducted in Section 6.

Note that in Theorem 5.1 the last row of the matrix $D_3$ is a row of zeros. This implies that the mixing time of $\mathbf{1}^T\mathbf{X}^d$ grows more rapidly than $O(d)$ as $d \to \infty$. In [8], heuristic arguments and extensive simulations show that the mixing time of $\mathbf{1}^T\mathbf{X}^d$ is in fact $O(d^2)$. Theorem 5.3 below gives a formal statement of this result. (The proof of Theorem 5.3 is similar to the proof of Theorem 5.1 and is, hence, omitted.)

For $d \geq 1$, let $\tilde{U}_t^d = \frac{1}{d-2}\sum_{i=3}^d X_{[d^2t],i}^d$.

Theorem 5.3. *Suppose that $0 < \rho < 1$ and that $c_d \to c$, as $d \to \infty$, for some $0 < c \leq 1$. Let $\mathbf{X}_0^\infty = (X_{0,1}^1, X_{0,2}^2, \ldots)$ be constructed as in the prelude to Theorem 5.1. Then, as $d \to \infty$,*

$$\tilde{U}^d \Rightarrow \tilde{U},$$

*where $\tilde{U}_0 \sim N(0,\rho)$ and $\tilde{U}$ satisfies the Langevin SDE*

$$d\tilde{U}_t = (h_{c,\rho}(l))^{1/2} dB_t + h_{c,\rho}(l)\left\{-\frac{1}{2\rho}\tilde{U}_t\right\} dt,$$

*where* $h_{c,\rho}(l) = 2cl^2\Phi(-\frac{l}{2}\sqrt{\frac{c}{1-\rho}})$, *as before.*

We now turn our attention to MALA-within-Gibbs for the exchangeable normal. So that now the proposal $\mathbf{Y}^d$ is given by

$$Y_i^d = x_i^d + \chi_i^d\left\{\sigma_d Z_i + \frac{\sigma_d^2}{2}\left(-\frac{1}{1-\rho}x_i^d - \theta_d \sum_{j=1}^{d} x_j^d\right)\right\},$$

where we take $\sigma_d^2 = l^2 d^{-1/3}$ with $l$ an arbitrary constant. Let $\mathbf{X}_0^\infty$ be constructed as outlined above for the RWM-within-Gibbs. Let $\{J_t\}$ be a Poisson process with rate $d^{1/3}$ and let $\Gamma^d = \{\Gamma_t^d\}_{t\geq 0}$ be the $d$-dimensional jump process defined by $\Gamma_t^d = \mathbf{X}_{J_t}^d$. Let $\mathbf{U}_t^d = (U_{t,1}^d, U_{t,2}^d, U_{t,3}^d)$ be such that $U_{t,1}^d = \Gamma_{t,1}^d$, $U_{t,2}^d = \Gamma_{t,2}^d$ and $U_{t,3}^d = \frac{1}{d-2}\sum_{i=3}^d \Gamma_{t,i}$.

THEOREM 5.4. *Suppose that* $0 < \rho < 1$ *and that* $c_d \to c$, *as* $d \to \infty$, *for some* $0 < c \leq 1$. *Let* $\mathbf{X}_0^\infty = (X_{0,1}^1, X_{0,2}^2, \ldots)$ *be constructed as in the prelude to Theorem 5.1. Let* $D_1$, $D_2$, $D_3$ *and* $\tilde{f}$ *be as defined in Theorem 5.1. Then, as* $d \to \infty$,

$$\mathbf{U}^d \Rightarrow \mathbf{U},$$

*where* $\mathbf{U}_0$ *is distributed according to* $\tilde{f}$ *and* $\mathbf{U}$ *satisfies the Langevin SDE*

$$d\mathbf{U}_t = (h_{c,\rho}(l))^{1/2} D_3\, d\mathbf{B}_t + h_{c,\rho}(l) D_3\{\tfrac{1}{2}\nabla(-\tfrac{1}{2}\mathbf{U}_t^T D_2 \mathbf{U}_t)\}\, dt,$$

*where*

$$h_{c,\rho}(l) = 2cl^2\Phi\left(-\frac{l^3}{8}\sqrt{\frac{c}{(1-\rho)^3}}\right)$$

*is the speed of the limiting diffusion.*

Note that if we define

$$\tilde{K}_d^2 = \mathbb{E}\left[\frac{1}{48}\left\{5\left(\frac{\partial^3}{\partial x_1^3} j(\mathbf{X}^d)\right)^2 - 3\left(\frac{\partial^3}{\partial x_1^2} j(\mathbf{X}^d)\right)^3\right\}\right]$$

and $\tilde{K}^2 = \frac{1}{16}(\frac{1}{1-\rho})^3$, then $\tilde{K}_d^2 \to \tilde{K}^2$ as $d \to \infty$ and $h_{c,\rho}(l) = 2cl^2\Phi(-\frac{l^3}{2}\sqrt{c}\tilde{K})$. Therefore, the speed of the limiting diffusion for exchangeable normal has the same form as that obtained for the IID product densities considered in Section 4.

As in (2.3), let $a_d^{c_d,\rho}(l)$ be the $\pi_d$-average acceptance rate of the above algorithm where $\mathbf{X}^d \sim N(\mathbf{0}, \Sigma_\rho^d)$, $\sigma_d = ld^{-1/6}$ and we update a proportion $c_d$ of the $d$ components in each iteration. Then we have the following corollary.

COROLLARY 5.5. *Let $c_d \to c$, as $d \to \infty$, for some $0 < c \leq 1$. Then, for $0 < \rho < 1$:*

(i) $\lim_{d\to\infty} a_d^{c_d,\rho}(l) = a^{c,\rho}(l) \stackrel{def}{=} 2\Phi(-\frac{l^3}{8}\sqrt{\frac{c}{(1-\rho)^3}})$.

(ii) *Let $\hat{l}$ be the unique value of $l$ which maximizes $h_{1,0}(l) = 2l^2\Phi(-\frac{l^3}{8})$ on $[0,\infty)$, and let $\hat{l}_{c,\rho}$ be the unique value of $l$ which maximizes $h_{c,\rho}(\hat{l})$ on $[0,\infty)$. Then $\hat{l}_{c,\rho} = \sqrt{1-\rho}c^{-1/6}l$ and $h_{c,\rho}(\hat{l}_{c,\rho}) = c^{2/3}(1-\rho)h_{1,0}(\hat{l})$.*

(iii) *For all $0 < c \leq 1$ and $0 < \rho < 1$, the optimal acceptance rate $a^{c,\rho}(\hat{l}_{c,\rho}) =$ 0.574 (to three decimal places).*

Note that Corollary 5.5(ii) states that the cost incurred by having $\sigma_{ij}^d = \rho$, $i \neq j$, rather than $\sigma_{ij}^d = 0$, $i \neq j$, is to slow down the speed of the limiting diffusion by a factor of $1-\rho$, for all $0 < c \leq 1$. Therefore, the dependence in the target distribution $\pi_d(\cdot)$ affects convergence of the MALA-within-Gibbs in the same way that it affects the RWM-within-Gibbs. The cost associated with updating only a proportion $c$ rather than all of the components is the same as that observed in Section 4. Furthermore, the optimal acceptance rate $a^{c,\rho}(\hat{l}_{c,\rho})$ is unaffected by the introduction of dependence.

From Theorem 5.4, we see that the mixing time of $\mathbf{1}^T\mathbf{X}^d$ is greater than $O(d^{1/3})$ as $d \to \infty$. In fact, the mixing time of $\mathbf{1}^T\mathbf{X}^d$ is in fact $O(d^{4/3})$. Let $\{J_t\}$ be a Poisson process with rate $d^{4/3}$ and for $d \geq 1$, let $\tilde{U}_t^d = \frac{1}{d-2}\sum_{i=3}^d X_{J_t,i}^d$.

THEOREM 5.6. *Suppose that $0 < \rho < 1$ and that $c_d \to c$, as $d \to \infty$, for some $0 < c \leq 1$. Let $\mathbf{X}_0^\infty = (X_{0,1}^1, X_{0,2}^2, \ldots)$ be constructed as in the prelude to Theorem 5.1. Then, as $d \to \infty$,*

$$\tilde{U}^d \Rightarrow \tilde{U},$$

*where $\tilde{U}_0 \sim N(0,\rho)$ and $\tilde{U}$ satisfies the Langevin SDE*

$$d\tilde{U}_t = (h_{c,\rho}(l))^{1/2}\,dB_t + h_{c,\rho}(l)\left\{-\frac{1}{2\rho}\tilde{U}_t\right\}dt,$$

*where $h_{c,\rho}(l) = 2cl^2\Phi(-\frac{l^3}{8}\sqrt{\frac{c}{(1-\rho)^3}})$, as before.*

The proofs of Theorems 5.4 and 5.6 are hybrids of those for the results of Section 4, and for Theorems 5.1 and 5.3 above, and are, hence, omitted.

**6. A simulation study.** The rotational symmetry of the Gaussian distribution effectively allows the dependence problem to be formulated as one of heterogeneity of scale. Other distributional forms exist for which this may be possible (e.g., the multivariate $t$-distribution), but it seems difficult to derive results for very general distributional families of target distribution without resorting to ideas such as this. Therefore, to support the conjecture that the conclusions of Sections 3–5 hold beyond the rigorous, theoretical results, we present the following simulation study. Furthermore, we demonstrate that the asymptotic results are achieved in relatively low dimensional ($d \geq 10$) situations.

Throughout the simulation study we measure speed/efficiency of the algorithm by considering first-order efficiency. That is, for a multidimensional Markov chain $\mathbf{X}$ with first component $X^1$, say, the first-order efficiency is defined to be $d\mathbb{E}[(X^1_{t+1} - X^1_t)^2]$ for RWM and $d^{1/3}\mathbb{E}[(X^1_{t+1} - X^1_t)^2]$ for MALA, where $\mathbf{X}_t$ is assumed to be stationary. For each of the target distributions and different choices of $c$ and $d$, we consider 50 different proposal variances, $\sigma^2_{d,c}$. For each choice of proposal variance $\sigma^2_{d,c}$, we started with $\mathbf{X}_0$ drawn from the target distribution. We then ran the algorithm for 100000 iterations. We estimate $\mathbb{E}[(X^1_{t+1} - X^1_t)^2]$ by $\frac{1}{100000}\sum_{i=1}^{100000}(X^1_i - X^1_{i-1})^2$ and the acceptance rate is estimated by $\frac{1}{100000}\sum_{i=1}^{100000}\mathbf{1}_{\{\mathbf{X}_i \neq \mathbf{X}_{i-1}\}}$. We then plot acceptance rate against $d\mathbb{E}[(X^1_{t+1} - X^1_t)^2]$ (first-order efficiency).

We begin by considering RWM-within-Gibbs. We shall consider three different target distributions $\pi_d \sim N(\mathbf{0}, \Sigma^d_\rho)$, $\pi_d \sim t_{50}(\mathbf{0}, \Sigma^d_\rho)$ and $\pi_d(\mathbf{x}^d) = \prod_{i=1}^d \frac{1}{2} \times \exp(-|x^d_i|)$ (double-sided exponential). Note that the distributions $t_{50}(\mathbf{0}, \Sigma^d_\rho)$ $(\rho > 0)$ and the double-sided exponential are not covered by the asymptotic results of Sections 3 and 5. For the $N(\mathbf{0}, \Sigma^d_\rho)$ and $t_{50}(\mathbf{0}, \Sigma^d_\rho)$, we plot acceptance rate against the normalized first-order efficiency, $\frac{d}{1-\rho}\mathbb{E}[(X^1_{t+1} - X^1_t)^2]$. The normalization is introduced to take account of dependence (see Corollary 5.2).

Figures 1 and 2 give a representative sample of the simulation study we conducted for a whole range of different values of $c$, $d$ and $\rho$. The results are as one would expect. In all cases the estimated optimal acceptance rate is approximately 0.234. As can be seen from Figures 1 and 2, the normalized first-order efficiency curves are virtually indistinguishable from one another for each choice of $c$, $d$ and $\rho$. Therefore, we have made no attempt to differentiate between the different efficiency curves.

(Note that the results in Figure 3 are a representative sample from a much larger simulation study.)

Figures 3 and 4 produce results in line with those expected from Sections 3 and 5. This demonstrates that the conclusions of Sections 3 and 5 do extend beyond those target distributions for which rigorous statements have been made.

We now turn our attention to MALA-within-Gibbs. We shall consider in our simulation study only target densities of the form $\pi_d \sim N(\mathbf{0}, \Sigma^d_\rho)$.

Simulations in Figures 5 and 6 show excellent agreement with Corollaries 4.3 and 5.5. Again, the results demonstrate the usefulness/relevance of the asymptotic results for even fairly small $d$.

**7. Discussion.** A rather surprising property of high-dimensional Metropolis and Langevin algorithms is the robustness of relative efficiency as a function of acceptance rate. In particular, the optimal acceptance rates 0.234 and 0.574 for Metropolis and Langevin, respectively, appear to be robust to many kinds of perturbation of the target density. A remarkable conclusion of this paper is this apparent robustness of relative efficiency, as a function of acceptance rate, seems to extend quite readily to updating schemes where only a fixed proportion of components are updated at once.

A further unexpected conclusion concerns the issue of optimization in $c$. Here, very clear cut statements appear to be available, with smaller-dimensional updates seeming to be optimal for the Metropolis algorithm (as seen from Theorem 3.1 and Corollary 3.2), whereas higher-dimensional up-dates are to be preferred (at least before computing time has been taken into consideration) for MALA schemes (see Theorem 4.2 and Corollary 4.3). The robustness of these conclusions to dependence in the target density is seen in the results of Section 5 and, supported by the simulation study in Section 6, seems contrary to the general intuition that "block updating" improves MCMC mixing (at least for the Metropolis results). However, our results show that this intuition is only correct for schemes where the multivariate update step utilies the structure of the target density (as, e.g., in the Gibbs sampler, or, to a lesser extent, MALA).

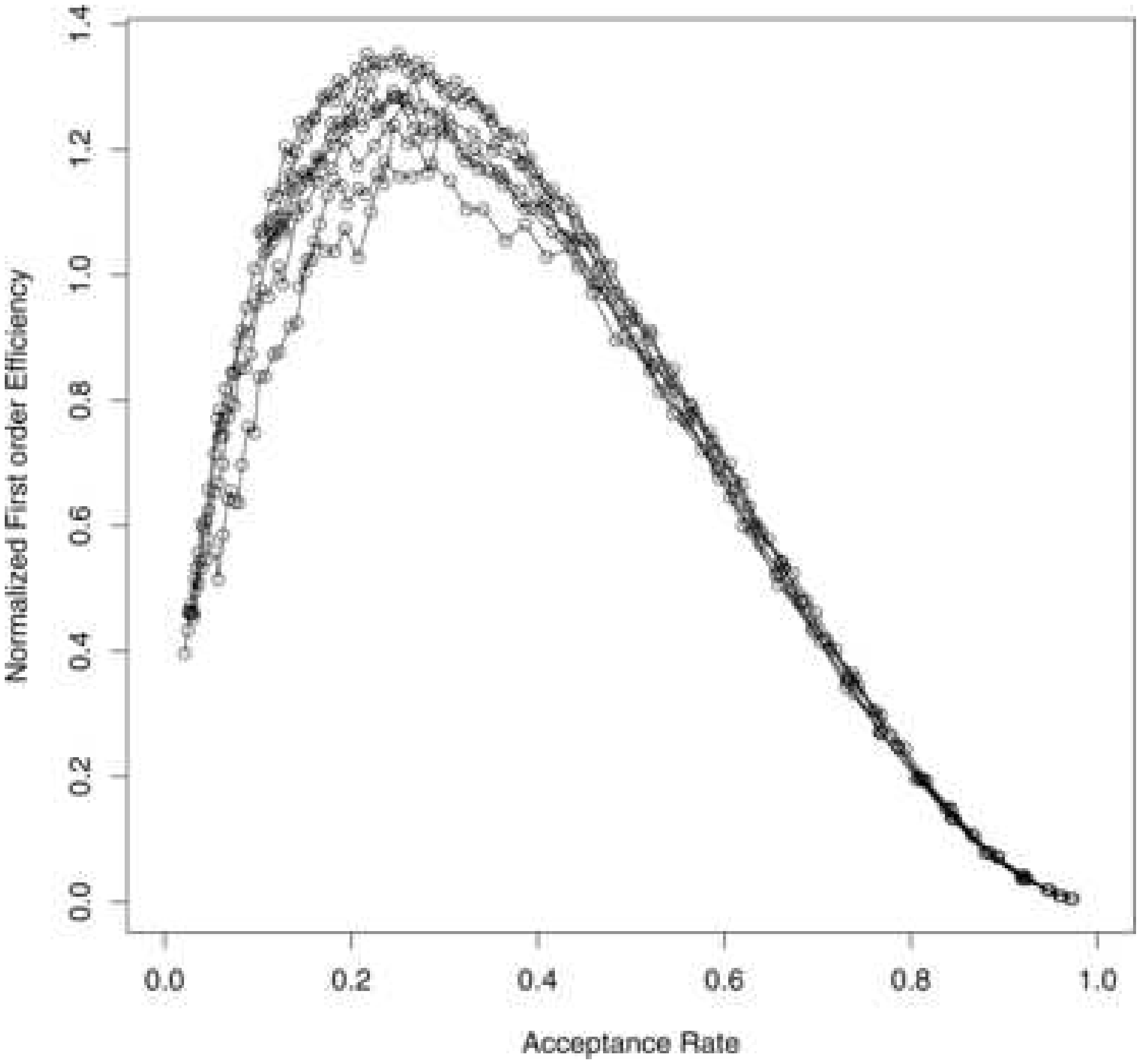

FIG. 1. *Normalized first-order efficiency of RWM-within-Gibbs,* $\frac{d}{1-\rho}\mathbb{E}[(X_{t+1}^1 - X_t^1)^2]$*, as a function of overall acceptance rates for each combination of* $(d = 20; c = 0.25, 0.5, 0.75, 1; \rho = 0, 0.5)$*, with* $\pi_d \sim N(\mathbf{0}, \Sigma_\rho^d)$.

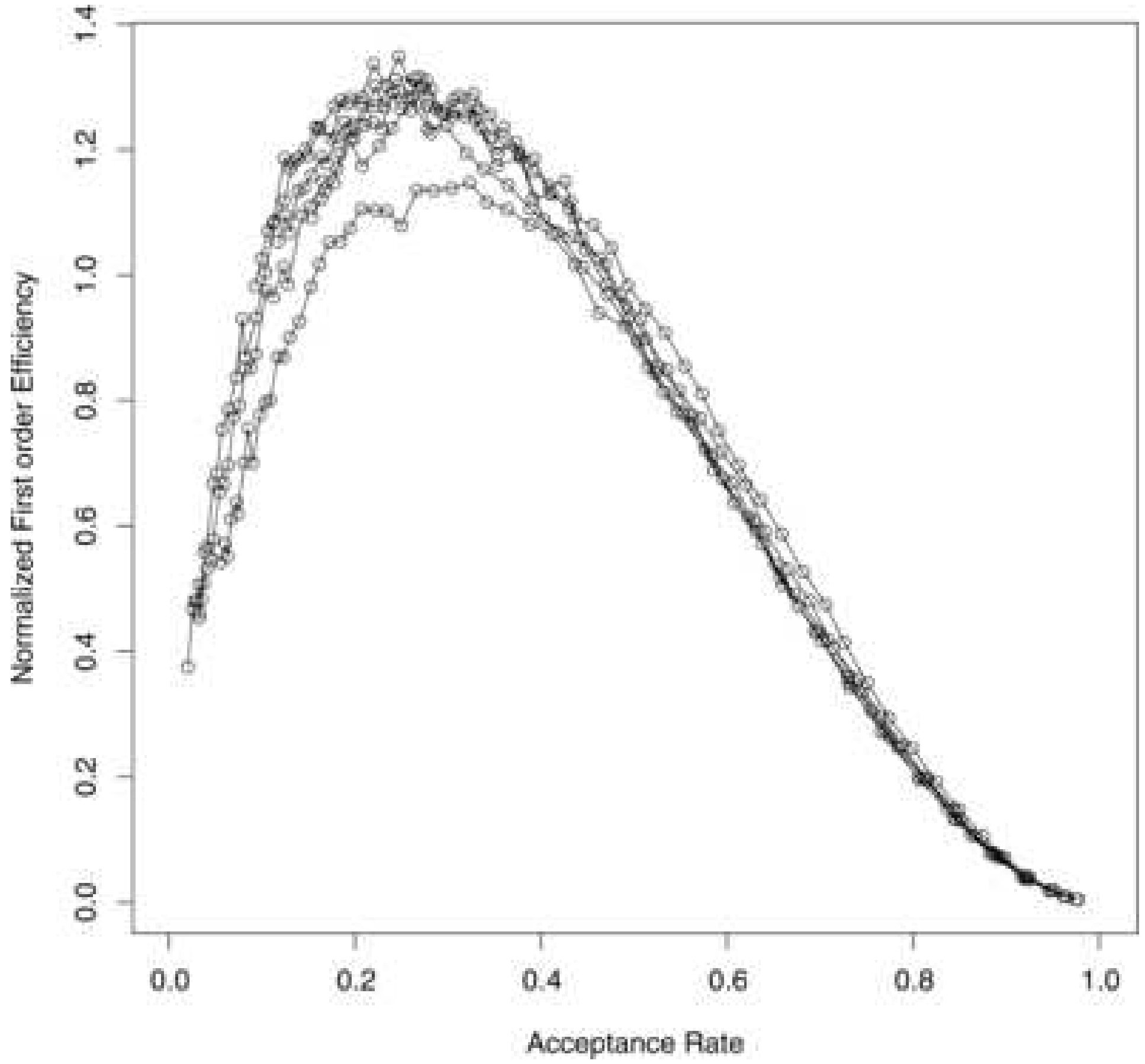

FIG. 2. *Normalized first-order efficiency of RWM-within-Gibbs,* $\frac{d}{1-\rho}\mathbb{E}[(X_{t+1}^1 - X_t^1)^2]$, *as a function of overall acceptance rates for each combination* ($c = 0.5$; $d = 10, 20, 50$; $\rho = 0, 0.5$), *with* $\pi_d \sim N(\mathbf{0}, \sigma_\rho^d)$.

We believe that these results should have quite fundamental implications for practical MCMC use, although, of course, they should be treated with care since they are only asymptotic. Our results have been shown in the simulation study to hold approximately in very low-dimensional problems—although the speed at which the infinite-dimensional limit is reached does vary in a complicated way, in particular, in $c$ and measures of dependence in the target density (such as $\rho$ in the exchangeable normal examples).

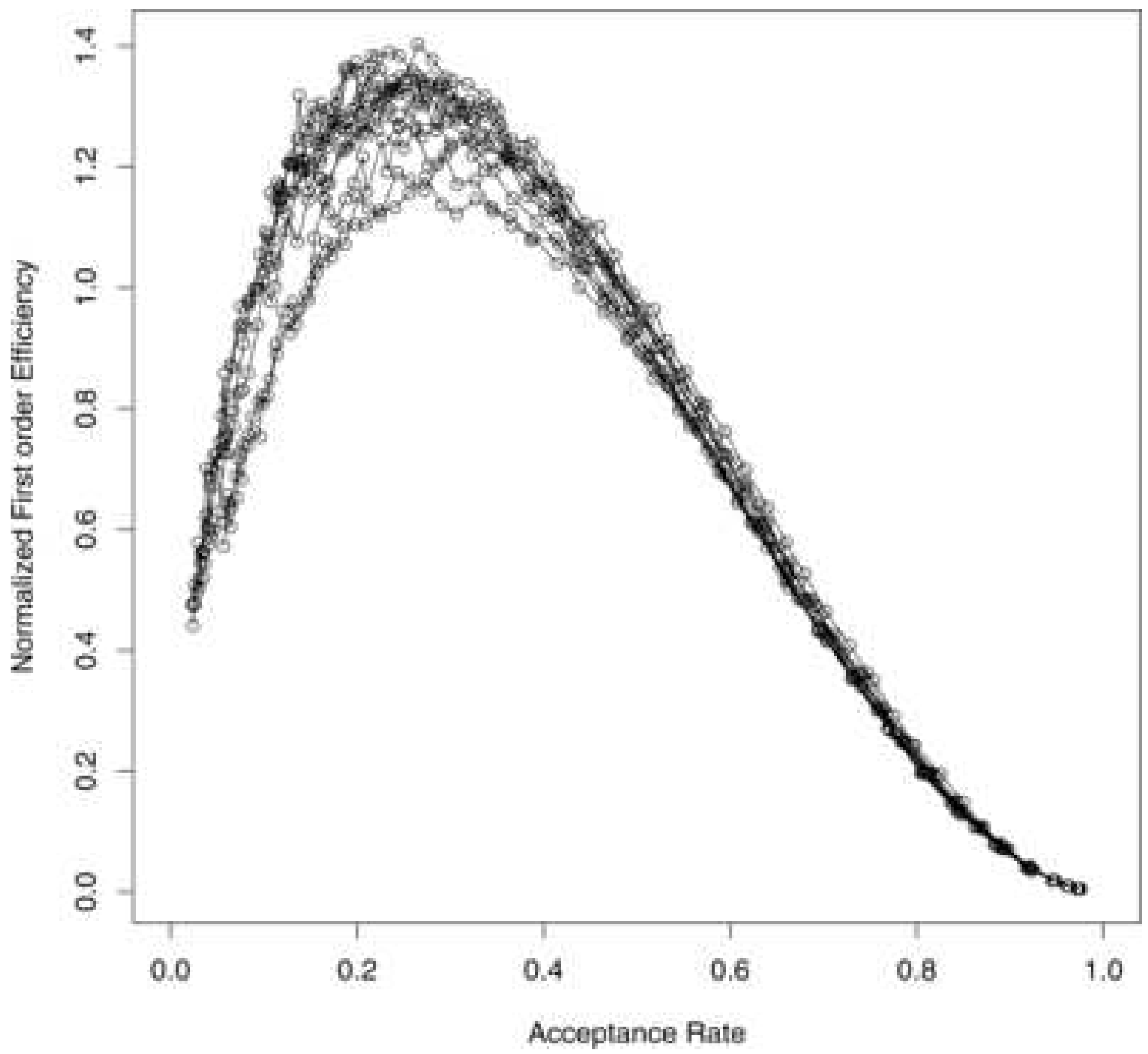

FIG. 3. *Normalized first-order efficiency of RWM-within-Gibbs,* $\frac{d}{1-\rho}\mathbb{E}[(X_{t+1}^1 - X_t^1)^2]$, *as a function of overall acceptance rates for each combination:* (i) $(d = 20;$ $c = 0.25, 0.5, 0.75, 1;$ $\rho = 0, 0.5)$ *and* (ii) $(c = 0.5;$ $d = 10, 20, 50;$ $\rho = 0, 0.5)$, *with* $\pi_d \sim t_{50}(\mathbf{0}, \Sigma_\rho^d)$.

The results for the exchangeable normal example show that certain functions can converge at different rates to others ($\bar{X}$ converging at rate $d^2$, while $X_i - \bar{X}$ converges at rate $d$), and this can cause serious practical problems for the MCMC practitioner. In particular, any one co-ordinate $X_i$ might converge rapidly, in a given time scale, to the wrong target density. Certainly, it would be extremely difficult to detect such problems empirically.

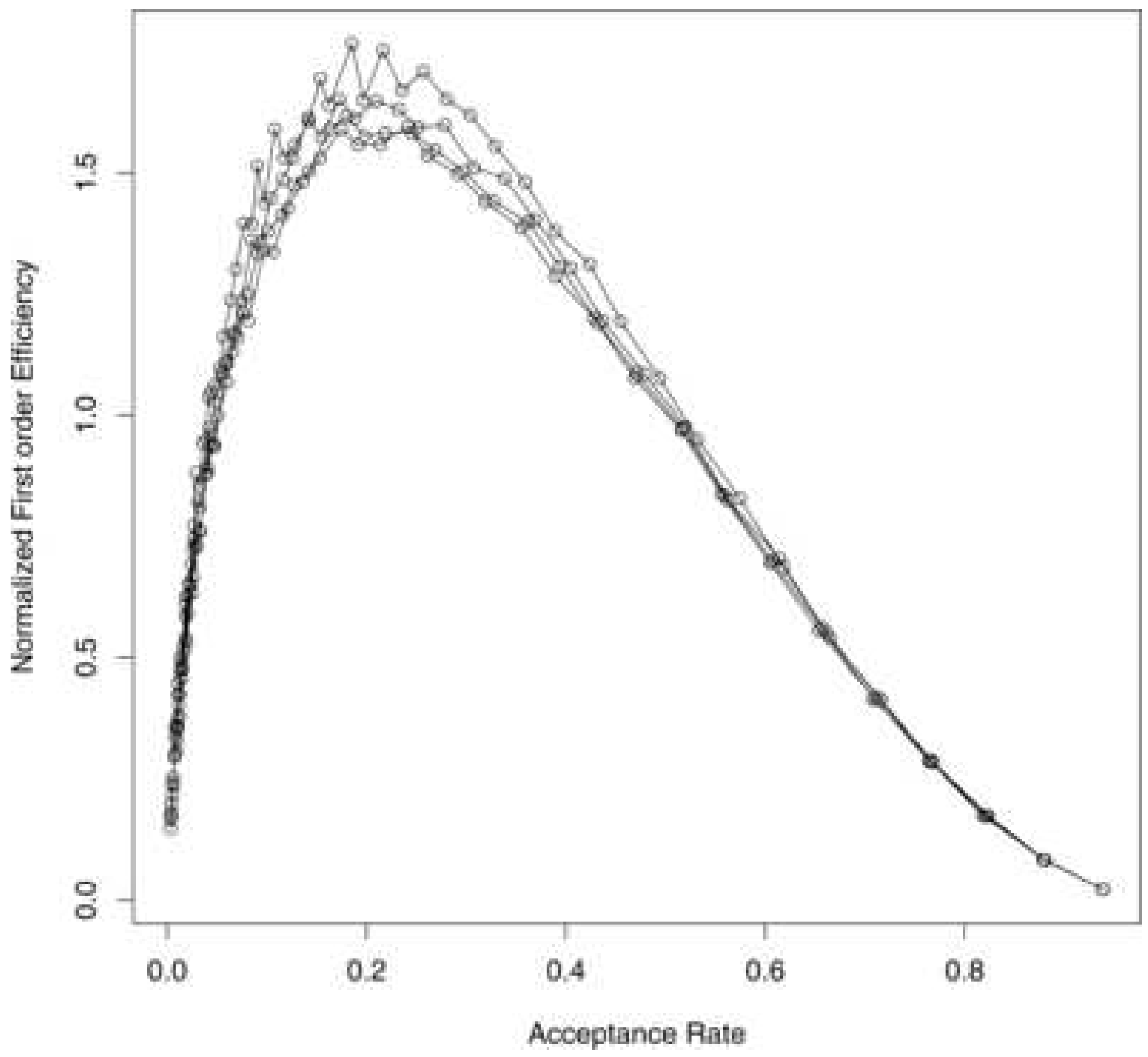

FIG. 4. *Normalized first-order efficiency of RWM-within-Gibbs,* $\frac{d}{1-\rho}\mathbb{E}[(X_{t+1}^1 - X_t^1)^2]$, *as a function of overall acceptance rates for each combination* $(d = 40; c = 0.25, 0.5, 0.75, 1)$, *with* $\pi_d(\mathbf{x}^d) = \prod_{i=1}^d \frac{1}{2}\exp(-|x_i^d|)$.

The results in this paper are given for Metropolis and MALA algorithms. However, the use of these two methods is, in some sense, illustrative, and other algorithms (such as, e.g., higher-order Langevin algorithms using, e.g., the Ozaki discretization [10]) are expected to yield similar conclusions.

## APPENDIX

**A.1. Proofs of Section 3.** Theorem 3.1 implies that the first component acts independently of all others as $d \to \infty$. Intuitively, this occurs because all other $(d-1)$ terms contribute expressions to the accept/reject ratio which turn out to obey SLLN and, thus, can be replaced by their deterministic limits. To make this idea rigorous, we need to define a set in $\mathbb{R}^d$ on which the first component is well approximated by the appropriate LLN limit.

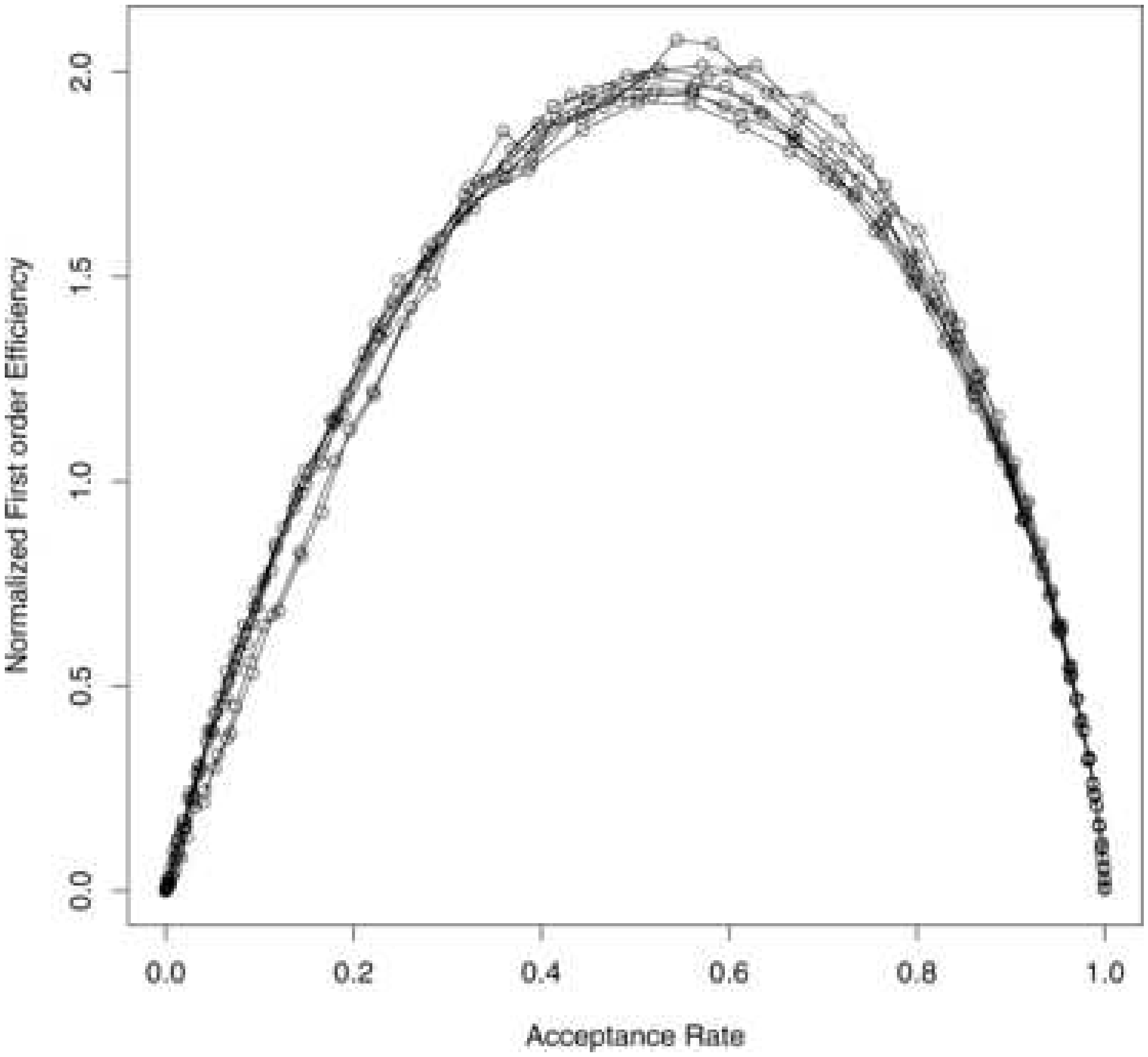

FIG. 5. *Normalized first-order efficiency of RWM-within-Gibbs,* $c^{-2/3}\frac{1}{1-\rho}d^{1/3}\times \mathbb{E}[(X^1_{t+1}-X^1_t)^2]$, *as a function of overall acceptance rates for each combination of* $(d=20;\ c=0.25, 0.5, 0.75, 1;\ \rho=0, 0.5)$, *with* $\pi_d \sim N(\mathbf{0}, \Sigma^d_\rho)$.

Motivated by this idea, we construct sets of tolerances around average values for quantities which will appear in the accept/reject ratio. Thus, we define the sequence of sets $\{F_d \subseteq \mathbb{R}^d, d > 1\}$ by

$$\begin{aligned} F_d = &\left\{\mathbf{x}^d; \left|\frac{1}{d-1}\sum_{i=2}^d g'(x_i^d)^2 - I\right| < d^{-1/8}\right\} \\ &\cap \left\{\mathbf{x}^d; \left|\frac{1}{d-1}\sum_{i=2}^d g''(x_i^d) + I\right| < d^{-1/8}\right\} \\ &\cap \left\{\mathbf{x}^d; \left|\frac{1}{(d-1)^2}\sum_{i=2}^d g'(x_i^d)^4\right| < d^{-1/8}\right\}, \end{aligned}$$

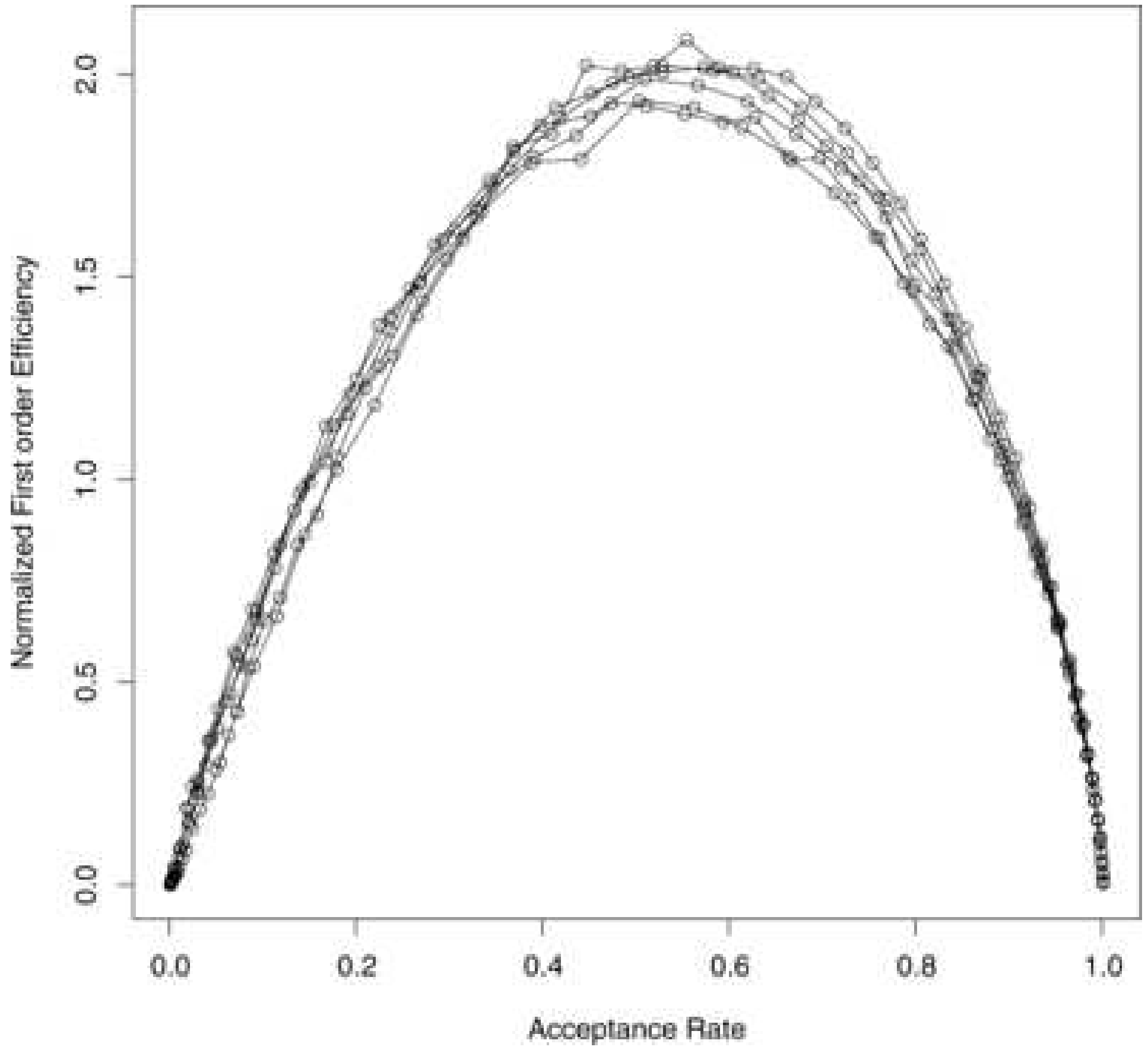

FIG. 6. *Normalized first-order efficiency of MALA-within-Gibbs,* $c^{-2/3}\frac{1}{1-\rho}d^{1/3} \times \mathbb{E}[(X^1_{t+1} - X^1_t)^2]$, *as a function of overall acceptance rates for each combination* ($c = 0.5$; $d = 10, 20, 50$; $\rho = 0, 0.5$), *with* $\pi_d \sim N(\mathbf{0}, \Sigma^d_\rho)$.

$$= F_{d,1} \cap F_{d,2} \cap F_{d,3}, \qquad \text{say,}$$

where $I$ is defined in Theorem 3.1. Let $\mathbf{x}^\infty = (x_1, x_2, \ldots)$ and for $d \geq 1$, let $\mathbf{x}^d = (x^d_1, x^d_2, \ldots, x^d_d)$, where, for $1 \leq i \leq d$, $x^d_i = x_i$. Thus, we shall use $x^d_1$ and $x_1$ interchangeably, as appropriate.

LEMMA A.1. *For* $k = 1, 2, 3$ *and* $t > 0$,

$$\mathbb{P}(\mathbf{U}^d_s \in F_{d,k}, 0 \leq s \leq t) \to 1 \qquad \text{as } d \to \infty \tag{A.1}$$

*and, hence,*

$$\mathbb{P}(\mathbf{U}^d_s \in F_d, 0 \leq s \leq t) \to 1 \qquad \text{as } d \to \infty.$$

PROOF. The cases $k=1$ and $k=2$ are proved in [6], Lemma 2.1. The case $k=3$ is proved similarly using Markov's inequality and (3.2). The lemma then follows. □

For any random variable $X$ and for any subset $A \subseteq \mathbb{R}$, let $\mathbb{E}^*[X] = \mathbb{E}[X|\chi_1^d = 1]$ and $\mathbb{P}^*(X \in A) = \mathbb{P}(X \in A|\chi_1^d = 1)$.

Let $G_d$ be the (discrete-time) generator of $\mathbf{X}^d$, and let $V \in C_c^\infty$ (the space of infinitely differentiable functions on compact support) be an arbitrary test function of the first component only. Thus,

$$\begin{aligned} G_d V(\mathbf{x}^d) &= d\mathbb{E}\left[(V(\mathbf{Y}^d) - V(\mathbf{x}^d))\left\{1 \wedge \frac{\pi_d(\mathbf{Y}^d)}{\pi_d(\mathbf{x}^d)}\right\}\right] \\ &= d\mathbb{P}(\chi_1^d = 1)\mathbb{E}^*\left[(V(\mathbf{Y}^d) - V(\mathbf{x}^d))\left\{1 \wedge \frac{\pi_d(\mathbf{Y}^d)}{\pi_d(\mathbf{x}^d)}\right\}\right], \end{aligned} \tag{A.2}$$

since $Y_1^d = x_1^d$ if $\chi_1^d = 0$.

The generator $G$ of the one-dimensional diffusion described in (3.4), for an arbitrary test function $V \in C_c^\infty$, is given by

$$GV(x_1) = 2cl^2\Phi\left(-\frac{l\sqrt{cI}}{2}\right)\left\{\frac{1}{2}g'(x_1)V'(x_1) + \frac{1}{2}V''(x_1)\right\}. \tag{A.3}$$

(Note that, under the conditions imposed in Theorem 3.1, $C_c^\infty$ forms a core for the full generator.) By Lemma A.1, we can restrict attention to $\mathbf{x}^d \in F_d$. The aim will therefore be to show that, for all $\mathbf{x}^d \in F_d$,

$$G_d V(\mathbf{x}^d) \to GV(x_1) \qquad \text{as } d \to \infty.$$

The proof of Theorem 3.1 will then be fairly straightforward.

Thus, we begin by giving a Taylor series approximation for $G_d V(\mathbf{x}^d)$ in Lemma A.3, for which we will require the following lemma.

LEMMA A.2. *For any $V \in C_c^\infty$ (the space of infinitely differentiable functions on compact support),*

$$\sup_{\mathbf{x}^d \in F_d} |d\mathbb{E}^*[(V(Y_1^d) - V(x_1^d))] - \tfrac{1}{2}l^2 V''(x_1)| \to 0 \qquad \textit{as } d \to \infty \tag{A.4}$$

*and*

$$\sup_{\mathbf{x}^d \in F_d} |\sigma_d d\mathbb{E}^*[Z_1(V(Y_1^d) - V(x_1^d))] - l^2 V'(x_1)| \to 0 \qquad \textit{as } d \to \infty, \tag{A.5}$$

*with $x_1 = x_1^d$.*

Proof. For $\chi_1^d = 1$,

$$Y_1^d - x_1^d = \sigma_d Z_1 = \frac{l}{\sqrt{d-1}} Z_1.$$

Thus, by Taylor's theorem,

$$\begin{aligned} V(Y_1^d) - V(x_1^d) &= V'(x_1^d)(\sigma_d Z_1) \\ &\quad + \tfrac{1}{2}V''(x_1^d)(\sigma_d Z_1)^2 + \tfrac{1}{6}V'''(W_1)(\sigma_d Z_1)^3 \end{aligned} \tag{A.6}$$

for some $W_1$ lying between $x_1^d$ and $Y_1^d$.

The lemma then follows by substituting (A.6) into the left-hand sides of (A.4) and (A.5). □

Lemma A.3. *Let*

$$\tilde{G}_d V(\mathbf{x}^d) = \tfrac{1}{2} c l^2 V''(x_1)\mathbb{E}^*[1 \wedge e^{B_d}] + c l^2 V'(x_1) g'(x_1)\mathbb{E}^*[e^{B_d}; B_d < 0],$$

*where* $B_d (= B_d(\mathbf{x}^d)) = \sum_{i=2}^d (g(Y_i^d) - g(x_i^d))$. *Then, we have that*

$$\sup_{\mathbf{x}^d \in F_d} |G_d V(\mathbf{x}^d) - \tilde{G}_d V(\mathbf{x}^d)| \to 0 \qquad \text{as } d \to \infty. \tag{A.7}$$

Proof. Decomposing $\mathbf{Y}^d$ into $(Y_1, \mathbf{Y}^{d-})$ and using independence gives

$$G_d V(\mathbf{x}^d) = d c_d \mathbb{E}^*_{Y_1^d}\left[(V(Y_1^d) - V(x_1^d))\mathbb{E}^*_{\mathbf{Y}^{d-}}\left[1 \wedge \prod_{i=1}^d \frac{f(Y_i^d)}{f(x_i^d)}\right]\right].$$

We shall begin by concentrating on the inner expectation, by recalling the following fact noted in [2]. Let $h$ be a twice differentiable function on $\mathbb{R}$, then the function $z \mapsto 1 \wedge e^{h(z)}$ is also twice differentiable, except at a countable number of points, with first derivative given Lebesgue almost everywhere by the function

$$\frac{d}{dz} 1 \wedge e^{h(z)} = \begin{cases} h'(z)e^{h(z)}, & \text{if } h(z) < 0, \\ 0, & \text{if } h(z) \geq 0. \end{cases}$$

Now take $h_d(z) (= h_d(z; \mathbf{x}^d)) = (g(x_1^d + \sigma_d z) - g(x_1^d)) + B_d$ and let

$$\gamma_d(z) = \mathbb{E}^*_{\mathbf{Y}^{d-}}\left[1 \wedge \prod_{i=1}^d \frac{f(Y_i^d)}{f(x_i^d)} \middle| Z_1 = z\right].$$

Thus, $\gamma_d(z) = \mathbb{E}^*_{\mathbf{Y}^{d-}}[1 \wedge e^{h_d(z)}]$, and so, for almost every $x_1^d \in \mathbb{R}$, there exists $W$ lying between 0 and $z$ such that

$$\begin{aligned} \gamma_d(z) &= \mathbb{E}^*_{\mathbf{Y}^{d-}}[1 \wedge e^{h_d(0)}] \\ &\quad + z\mathbb{E}^*_{\mathbf{Y}^{d-}}[\sigma_d g'(x_1^d) e^{h_d(0)}; h_d(0) < 0] \\ &\quad + \frac{z^2}{2}\mathbb{E}^*_{\mathbf{Y}^{d-}}[\sigma_d^2 (g''(x_1^d + \sigma_d W) + g'(x_1^d + \sigma_d W)^2) e^{h_d(W)}; h_d(W) < 0]. \end{aligned} \tag{A.8}$$

The key results to note are that $h_d(0) = B_d$ and that, conditional upon $\chi_1^d = 1$, $Y_1^d$ and $\mathbf{Y}^{d-}$ are independent. Therefore,

$$\begin{aligned}
&G_d V(\mathbf{x}^d)\\
&\quad = dc_d \mathbb{E}^*_{Y_1^d}\Big[(V(Y_1^d) - V(x_1^d))\\
&\qquad\qquad \times \Big\{\mathbb{E}^*_{\mathbf{Y}^{d-}}[1 \wedge e^{h_d(0)}]\\
&\qquad\qquad\quad + Z_1 \mathbb{E}^*_{\mathbf{Y}^{d-}}[\sigma_d g'(x_1^d) e^{h_d(0)}; h_d(0) < 0]\\
&\qquad\qquad\quad + \frac{Z_1^2}{2}\mathbb{E}^*_{\mathbf{Y}^{d-}}[\sigma_d^2(g''(x_1^d + \sigma_d W)\\
&\qquad\qquad\qquad\qquad + g'(x_1^d + \sigma_d W)^2) e^{h_d(W)}; h_d(W) < 0]\Big\}\Big]\\
&\quad = dc_d \mathbb{E}^*[(V(Y_1^d) - V(x_1))]\mathbb{E}^*[1 \wedge e^{B_d}]\\
&\qquad + g'(x_1) dc_d \sigma_d \mathbb{E}^*[(V(Y_1^d) - V(x_1^d)) Z_1]\mathbb{E}^*[e^{B_d}; B_d < 0] \qquad \text{(A.9)}\\
&\qquad + dc_d \mathbb{E}^*_{Y_1^d}\Big[(V(Y_1^d) - V(x_1^d))\frac{Z_1^2}{2}\\
&\qquad\qquad \times \mathbb{E}^*_{\mathbf{Y}^{d-}}[\sigma_d^2(g''(x_1^d + \sigma_d W)\\
&\qquad\qquad\qquad + g'(x_1^d + \sigma_d W)^2) e^{h_d(W)}; h_d(W) < 0]\Big]\\
&\quad = \hat{G}_d V(\mathbf{x}^d) + D_d(\mathbf{x}^d; Z_1; W), \qquad \text{say.}
\end{aligned}$$

Since $\mathbb{E}^*[1 \wedge e^{B_d}], \mathbb{E}^*[e^{B_d}; B_d < 0] \le 1$ and $x_1 = x_1^d$, it follows from Lemma A.2 that

$$\sup_{\mathbf{x^d} \in F_d} |\hat{G}_d V(\mathbf{x}^d) - \tilde{G}_d V(\mathbf{x}^d)| \to 0 \qquad \text{as } d \to \infty.$$

Thus, to prove the lemma, it is sufficient to show that, for all $\mathbf{x}^d \in F_d$, $D_d(\mathbf{x}^d; Z_1; W)$ converges to 0, as $d \to \infty$.

By Taylor's theorem, we have that

$$\left|(V(Y_1^d) - V(x_1^d))\frac{Z_1^2}{2}\right| \le \sup_{a_1 \in \mathbb{R}} |V'(a_1)| \frac{\sigma_d}{2} |Z_1^3|$$

and

$$\begin{aligned}
|g'(x_1^d + \sigma_d W)| &\le |g'(x_1^d)| + \sigma_d |W| \sup_{a_2 \in \mathbb{R}} |g''(a_2)|\\
&\le |g'(x_1^d)| + \sigma_d |Z_1| \sup_{a_2 \in \mathbb{R}} |g''(a_2)|.
\end{aligned}$$

Since $V'$ and $g''$ are bounded functions, it follows that, for all $\mathbf{x}^d \in F_d$,

$$D_d(\mathbf{x}^d; Z_1; W) \le d c_d \{\tfrac{3}{2} K \sigma_d^3 (K + |g'(x_1^d)| + \sigma_d K)\} \to 0 \qquad \text{as } d \to \infty,$$

for some $K > 0$, and the lemma is proved. $\square$

Lemma A.3 states that, for all $\mathbf{x}^d \in F_d$, the generator $G_d$ can be approximated by the generator $\tilde{G}_d$ which resembles the limiting generator $G$. Thus, we now need to consider for all $\mathbf{x}^d \in F_d$, $\mathbb{E}^*[1 \wedge e^{B_d}]$ and $\mathbb{E}^*[e^{B_d}; B_d < 0]$. The aim is to approximate $B_d$ by a more convenient quantity $A_d$ (to be defined in Lemma A.6) and, hence, show that

$$\mathbb{E}^*[1 \wedge e^{B_d}] \to 2\Phi\left(-\frac{l\sqrt{cI}}{2}\right)$$

and

$$\mathbb{E}^*[e^{B_d}; B_d < 0] \to \Phi\left(-\frac{l\sqrt{cI}}{2}\right) \qquad \text{as } d \to \infty.$$

This will be done in the following lemmas.

Lemma A.4. *Let* $\lambda_d (= \lambda_d(\mathbf{x}^d)) = \frac{1}{d-1} \sum_{i=2}^d \chi_i^d g'(x_i^d)^2$. *For any* $\varepsilon > 0$,

$$\sup_{\mathbf{x}^d \in F_d} \mathbb{P}^*(|\lambda_d - cI| > \varepsilon) \to 0 \qquad \textit{as } d \to \infty.$$

Proof. Let $R_d (= R_d(\mathbf{x}^d)) = \frac{1}{d-1} \sum_{i=2}^d g'(x_i^d)^2$. Then, for $\mathbf{x}^d \in F_d$,

$$|\lambda_d - cI| \le |\lambda_d - \mathbb{E}^*[\lambda_d]| + |\mathbb{E}^*[\lambda_d] - cR_d| + |cR_d - cI|.$$

Note that $\mathbb{E}^*[\lambda_d] = \frac{c_d d - 1}{d-1} R_d$, and so, by Lemma A.1, we have that

$$|\mathbb{E}^*[\lambda_d] - cR_d| + |cR_d - cI| \to 0 \qquad \text{as } d \to \infty.$$

Therefore, to prove the lemma, it suffices to show that, for any $\varepsilon > 0$, $\mathbb{P}^*(|\lambda_d - \mathbb{E}^*[\lambda_d]| > \varepsilon) \to 0$ as $d \to \infty$. Note that

$$\lambda_d^2 = \frac{1}{(d-1)^2} \sum_{i=2}^d \sum_{j=2}^d \chi_i^d \chi_j^d g'(x_i^d)^2 g'(x_j^d)^2,$$

and so,

$$\mathbb{E}^*[\lambda_d^2] = \frac{1}{(d-1)^2} \left\{ \frac{c_d d - 1}{d-1} \sum_{i=2}^d g'(x_i^d)^4 \right.$$
$$\left. + \frac{(c_d d - 1)}{(d-1)} \frac{(c_d d - 2)}{(d-2)} \sum_{i=2}^d \sum_{j \ne i} g'(x_i^d)^2 g'(x_j^d)^2 \right\}$$

$$= \frac{(c_d d-1)(c_d d-2)}{(d-1)(d-2)} R_d^2 + \frac{(c_d d-1)(1-c_d)d}{(d-1)(d-2)} \left\{ \frac{1}{(d-1)^2} \sum_{i=2}^{d} g'(x_i^d)^4 \right\}.$$

Then since $\sup_{\mathbf{x}^d \in F_d} |\frac{1}{(d-1)^2} \sum_{i=2}^d g'(x_i^d)^4| \to 0$ and $c_d \to c$ as $d \to \infty$, it follows that, for all $\mathbf{x}^d \in F_d$, $\mathbb{E}^*[(\lambda_d - \mathbb{E}^*[\lambda_d])^2] \to 0$ as $d \to \infty$ and, hence, by Chebyshev's inequality,

$$\sup_{\mathbf{x}^d \in F_d} \mathbb{P}^*(|\lambda_d - \mathbb{E}^*[\lambda_d]| > \varepsilon) \to 0 \qquad \text{as } d \to \infty,$$

as required. □

LEMMA A.5. *Let*

$$W_d (= W_d(\mathbf{x}^d)) = \sum_{i=2}^{d} \left\{ \frac{1}{2} g''(x_i^d)(Y_i^d - x_i^d)^2 + \frac{cl^2}{2(d-1)} g'(x_i^d)^2 \right\},$$

*and* $c_d \to c$ *as* $d \to \infty$. *Then, recalling that* $\sigma_d = l/\sqrt{d}$,

$$\sup_{\mathbf{x}^d \in F_d} \mathbb{E}^*[|W_d|] \to 0 \qquad \text{as } d \to \infty.$$

PROOF. First, note that $\mathbb{E}^*[|W_d|]^2 \le \mathbb{E}^*[W_d^2]$.

Then, by direct calculations,

$$\begin{aligned}
\mathbb{E}^*[W_d^2] &= \sum_{i=2}^{d} \sum_{j=2}^{d} \mathbb{E}^* \left[ \left\{ \frac{1}{2} g''(x_i^d)(Y_i^d - x_i^d)^2 + \frac{cl^2}{2(d-1)} g'(x_i^d)^2 \right\} \right. \\
&\qquad \left. \times \left\{ \frac{1}{2} g''(x_j^d)(Y_j^d - x_j^d)^2 + \frac{cl^2}{2(d-1)} g'(x_j^d)^2 \right\} \right] \\
&= \left( \sum_{i=2}^{d} \frac{1}{4} g''(x_i^d)^2 \left\{ 3 \frac{c_d d - 1}{d-1} \sigma_d^4 - \frac{(c_d d-1)(c_d d-2)}{(d-1)(d-2)} \sigma_d^4 \right\} \right) \\
&\quad + \left( \sum_{i=2}^{d} \sum_{j=2}^{d} \left\{ \frac{1}{4} g''(x_i^d) g''(x_j^d) \frac{(c_d d-1)(c_d d-2)}{(d-1)(d-2)} \sigma_d^4 \right. \right. \\
&\qquad + \frac{cl^2}{2(d-1)} g''(x_i^d) g'(x_j^d)^2 \frac{c_d d-1}{d-1} \sigma_d^2 \\
&\qquad \left. \left. + \frac{c^2 l^4}{4(d-1)^2} g'(x_i^d)^2 g'(x_j^d)^2 \right\} \right) \\
&= W_{d,1} + W_{d,2}, \qquad \text{say.}
\end{aligned}$$

Let $W_{d,3}(= W_{d,3}(\mathbf{x}^d)) = \{\frac{cl^2}{2(d-1)} \sum_{i=2}^d (g''(x_i^d) + g'(x_i^d)^2)\}^2$, and since $c_d \to c$ as $d \to \infty$, we have that

$$\sup_{\mathbf{x}^d \in F_d} |W_{d,2} - W_{d,3}| \to 0 \qquad \text{as } d \to \infty.$$

However, by definition, $\sup_{\mathbf{x}^d \in F_d} |W_{d,3}| \to 0$ and since $g''$ is bounded, $\sup_{\mathbf{x}^d \in F_d} |W_{d,1}| \to 0$ as $d \to \infty$. The lemma follows immediately. □

Lemma A.6. *Let $A_d(= A_d(\mathbf{x}^d)) = \sum_{i=2}^d \{g'(x_i^d)(Y_i^d - x_i^d) - \frac{cl^2}{2(d-1)} g'(x_i^d)^2\}$. Then,*

$$\text{(A.10)} \qquad \sup_{\mathbf{x}^d \in F_d} |\mathbb{E}^*[1 \wedge e^{A_d}] - \mathbb{E}^*[1 \wedge e^{B_d}]| \to 0 \qquad \text{as } d \to \infty$$

*and*

$$\text{(A.11)} \qquad \sup_{\mathbf{x}^d \in F_d} |\mathbb{E}^*[e^{A_d}; A_d < 0] - \mathbb{E}^*[e^{B_d}; B_d < 0]| \to 0 \qquad \text{as } d \to \infty.$$

Proof. Note that

$$\begin{aligned} B_d &= \sum_{i=2}^d (g(Y_i^d) - g(x_i^d)) \\ &= \sum_{i=2}^d \{g'(x_i^d)(Y_i^d - x_i^d) + \tfrac{1}{2} g''(x_i^d)(Y_i^d - x_i^d)^2 + \tfrac{1}{6} g'''(\alpha_i^d)(Y_i^d - x_i^d)^3\}, \end{aligned}$$

for some $\alpha_i^d$ lying between $x_i^d$ and $Y_i^d$. Therefore, by [6], Proposition 2.2,

$$\begin{aligned} &|\mathbb{E}^*[1 \wedge e^{A_d}] - \mathbb{E}^*[1 \wedge e^{B_d}]| \\ \text{(A.12)} \qquad &\le \mathbb{E}^*[|W_d|] + \sup_{a \in \mathbb{R}} |g'''(a)| \frac{d-1}{6} \mathbb{E}^*[|Y_2^d - x_2^d|^3], \\ &= \mathbb{E}^*[|W_d|] + \sup_{a \in \mathbb{R}} |g'''(a)| \frac{d-1}{6} \frac{c_d d - 1}{d-1} \sigma_d^3 \mathbb{E}[|Z_1|^3]. \end{aligned}$$

Now let $\varphi_d = \sup_{\mathbf{x}^d \in F_d} \{\mathbb{E}^*[|W_d|] + \sup_{a \in \mathbb{R}} |g'''(a)| \frac{c_d d - 1}{6} \sigma_d^3 \mathbb{E}[|Z_1|^3]\}$, where $W_d$ is defined in Lemma A.5. Then, since $g'''$ is a bounded function, it follows from Lemma A.5 that $\varphi_d \to 0$ as $d \to \infty$ and so (A.10) is proved.

Let $J_d(= J_d(\mathbf{x}^d)) = (e^{A_d}; A_d < 0) - (e^{B_d}; B_d < 0)$ and let $\delta_d = \sqrt{\varphi_d}$. Then we proceed by showing that

$$\sup_{\mathbf{x}^d \in F_d} \mathbb{P}^*(|J_d| > \delta_d) \to 0 \qquad \text{as } d \to \infty.$$

Note that, if $A_d, B_d > 0$, then

$$|J_d| = 0 \le |A_d - B_d|$$

and if $A_d, B_d < 0$, then

$$|J_d| = |\exp(A_d) - \exp(B_d)| \leq |A_d - B_d|.$$

Therefore, it follows that

$$\mathbb{P}^*(|J_d| > \delta_d) \leq \mathbb{P}^*(-\delta_d < A_d < \delta_d) + \mathbb{P}^*(|A_d - B_d| \geq \delta_d). \tag{A.13}$$

By Markov's inequality,

$$\begin{aligned} &\mathbb{P}^*(|A_d - B_d| \geq \delta_d) \\ &\quad \leq \frac{1}{\delta_d}\mathbb{E}^*[|A_d - B_d|] \\ &\quad \leq \frac{1}{\delta_d}\left\{\mathbb{E}^*[|W_d|] + \sup_{a\in\mathbb{R}}|g'''(a)|\frac{d-1}{6}\mathbb{E}^*[|Y_2 - x_2|^3]\right\} \\ &\quad \leq \sqrt{\varphi_d}, \end{aligned} \tag{A.14}$$

and so, $\mathbb{P}^*(|A_d - B_d| > \delta_d) \to 0$ as $d \to \infty$, uniformly for $\mathbf{x}^d \in F_d$.

Fix $\mathbf{x}^d \in F_d$, then for any $\varepsilon > 0$, by Lemma A.4,

$$\mathbb{P}^*\left(\left|\frac{1}{l\sqrt{\lambda_d}}\left(\pm\delta_d + \frac{cl^2}{2}R_d\right) - \frac{l\sqrt{cI}}{2}\right| > \varepsilon\right) \to 0 \qquad \text{as } d \to \infty. \tag{A.15}$$

Hence,

$$\mathbb{E}^*\left[\Phi\left(\frac{1}{l\sqrt{\lambda_d}}\left(\pm\delta_d + \frac{cl^2}{2}R_d\right)\right)\right] \to \Phi\left(\frac{l\sqrt{cI}}{2}\right) \qquad \text{as } d \to \infty.$$

Thus,

$$\sup_{\mathbf{x}^d\in F_d} \mathbb{P}^*(-\delta_d < A_d < \delta_d) \to 0 \qquad \text{as } d \to \infty. \tag{A.16}$$

Therefore, by (A.13)–(A.16), $\sup_{\mathbf{x}^d\in F_d} \mathbb{P}^*(|J_d| > \delta_d) \to 0$ as $d \to \infty$. Then since $|J_d| \leq 1$, it follows that $\sup_{\mathbf{x}^d\in F_d} \mathbb{E}^*[J_d] \to 0$ as $d \to \infty$ and so (A.11) is proved. □

Lemma A.7.

$$\sup_{\mathbf{x}^d\in F_d}\left|\mathbb{E}^*[1 \wedge e^{A_d}] - 2\Phi\left(-\frac{l\sqrt{cI}}{2}\right)\right| \to 0 \qquad \text{as } d \to \infty \tag{A.17}$$

*and*

$$\sup_{\mathbf{x}^d\in F_d}\left|\mathbb{E}^*[e^{A_d}; A_d < 0] - \Phi\left(-\frac{l\sqrt{cI}}{2}\right)\right| \to 0 \qquad \text{as } d \to \infty. \tag{A.18}$$

Proof. Since $A_d \sim N(-\frac{cl^2}{2}R_d, l^2\lambda_d)$, it follows by [6], Proposition 2.4, that

$$\begin{aligned}\mathbb{E}^*[1 \wedge e^{A_d}] = \mathbb{E}^*\biggl[\Phi\biggl(-\frac{clR_d}{2\sqrt{\lambda_d}}\biggr) \\ + \exp\biggl(-\frac{l^2}{2}(cR_d - \lambda_d)\biggr)\Phi\biggl(-l\sqrt{\lambda_d} + \frac{clR_d}{2\sqrt{\lambda_d}}\biggr)\biggr].\end{aligned} \tag{A.19}$$

Since for any $\mathbf{x}^d \in F_d$ and $\varepsilon > 0$, $\mathbb{P}^*(|R_d - I| > \varepsilon) \to 0$ and $\mathbb{P}^*(|\lambda_d - cI| > \varepsilon) \to 0$ as $d \to \infty$, (A.17) follows from (A.19).

(A.18) is proved similarly. □

We are now in a position to show that, for all $\mathbf{x}^d \in F_d$, the generator $G_d$ converges to the generator $G$ as $d \to \infty$.

Theorem A.8. *For $V \in C_c^\infty$,*

$$\sup_{\mathbf{x}^d \in F_d} |G_d V(\mathbf{x}^d) - GV(x_1)| \to 0 \qquad \textit{as } d \to \infty.$$

Proof. By Lemma A.3,

$$\sup_{\mathbf{x}^d \in F_d} |G_d V(\mathbf{x}^d) - \tilde{G}V(\mathbf{x}^d)| \to 0 \qquad \text{as } d \to \infty,$$

and by Lemmas A.6 and A.7,

$$\sup_{\mathbf{x}^d \in F_d} |\tilde{G}_d V(\mathbf{x}^d) - GV(x_1)| \to 0 \qquad \text{as } d \to \infty.$$

Thus, the theorem is proved. □

Proof of Theorem 3.1. The proof is similar to that of [6]. From Lemmas A.1, A.4 and Theorem A.8, we have uniform convergence of $G_dV$ to $GV$ for vectors contained in a set of $\pi$ measure arbitrarily close to 1. Since $C_c^\infty$ separates points (see [4], page 113), the result will follow by [4], Chapter 4, Corollary 8.7 if we can demonstrate the compact containment condition, which in our case follows from the following statement. For all $\varepsilon > 0$, and all real valued $U_0^d = X_{0,1}^d$, we can find $K > 0$ sufficiently large with

$$\mathbb{P}(U_t^d \notin (-K, K),\ 0 \le t \le 1) \le \varepsilon,$$

for all $d$. We appeal directly to the explicit form of the Metropolis transitions and assume that the Lipshitz constant for $g$ is termed $b$. Thus, the following estimates are easy to derive by just noting that squared jumping distances are bounded above by that attained by ignoring rejections. Moreover, these estimates are uniform over all $\mathbf{X}_n^d$:

$$-b\sigma_d^2 e^{b^2\sigma_d^2/2} \le \mathbb{E}[X_{n+1,1}^d - X_{n,1}^d | \mathbf{X}_n^d] \le b\sigma_d^2 e^{b^2\sigma_d^2/2}$$

and

$$\mathbb{E}[(X_{n+1,1}^d - X_{n,1}^d)^2 | \mathbf{X}_n^d] \le \mathbb{E}[(Y_{n+1,1}^d - X_{n,1}^d)^2 | \mathbf{X}_n^d] = \sigma_d^2.$$

Thus, setting $V_n = X_{n,1}^d + nb\sigma_d^2 e^{b^2\sigma_d^2/2}$, $\{V_n, 0 \le n \le [d]\}$ is submartingale with

$$\mathbb{E}[V_{[d]}^2] \le d\sigma_d^2 + (db\sigma_d^2 e^{b^2\sigma_d^2/2})^2. \tag{A.20}$$

Since $\sigma_d^2 = \ell^2/d$, the right-hand side of (A.20) is uniformly bounded in $d$ so that the upper bound result follows by Doob's inequality. The lower bound follows similarly by considering the supermartingale $X_{n,1}^d - nb\sigma_d^2 e^{b^2\sigma_d^2/2}$. □

**A.2. Proofs of Section 4.** The proofs of Theorems 4.1 and 4.2 are similar to the proofs of Theorems 1 and 2 in [7], respectively. The only complication in the proofs is that we are updating a random set of components at each iteration in the MALA algorithm.

Let $\mathbf{x}^\infty = (x_1, x_2, \dots)$ and for $d \ge 1$, let $\mathbf{x}^d = (x_1^d, x_2^d, \dots, x_d^d)$, where, for $1 \le i \le d$, $x_i^d = x_i$. Thus, we shall again use $x_1^d$ and $x_1$ interchangeably as appropriate. Let $G_d$ be the (discrete-time) generator of $\mathbf{X}^d$ and let $V \in C_c^\infty$ be an arbitrary test function of the first component only. Thus,

$$\begin{aligned} G_d V(\mathbf{x}^d) &= d^{1/3}\mathbb{E}\left[(V(\mathbf{Y}^d) - V(\mathbf{x}^d))\left\{1 \wedge \frac{\pi_d(\mathbf{Y}^d)}{\pi_d(\mathbf{x}^d)}\right\}\right] \\ &= d^{1/3}\mathbb{P}(\chi_1^d = 1)\mathbb{E}^*\left[(V(\mathbf{Y}^d) - V(\mathbf{x}^d))\left\{1 \wedge \frac{\pi_d(\mathbf{Y}^d)}{\pi_d(\mathbf{x}^d)}\right\}\right], \end{aligned} \tag{A.21}$$

where $\mathbb{E}^*$ is defined after Lemma A.1 (cf. Section A.1 after Lemma A.1).

The generator $G$ of the one-dimensional diffusion described in Theorem 4.2, for an arbitrary test function, $V$, is given by

$$\begin{aligned} GV(x_1) &= 2cl^2\Phi\left(-\frac{\sqrt{c}l^3K}{2}\right)\left\{\frac{1}{2}g'(x_1)V'(x_1) + \frac{1}{2}V''(x_1)\right\} \\ &= h_c(l)\left\{\frac{1}{2}g'(x_1)V'(x_1) + \frac{1}{2}V''(x_1)\right\}, \end{aligned} \tag{A.22}$$

where $K$ and $h_c(l)$ are defined in Section 5.

The aim thus, as in Section A.1, is to find a sequence of sets $\{F_d \subseteq \mathbb{R}^d\}$ such that, for all $t > 0$,

$$\mathbb{P}(\Gamma_s^d \in F_d, \text{ for all } 0 \le s \le t) \to 1 \qquad \text{as } d \to \infty,$$

and, for $V \in C_c^\infty$,

$$\sup_{\mathbf{x}^d \in F_d} |G_d V(\mathbf{x}^d) - GV(x_1)| \to 0 \qquad \text{as } d \to \infty.$$

The proofs of Theorem 4.1 and 4.2 are then straightforward.

The first step is therefore to construct the sets $\{F_d \subseteq \mathbb{R}^d\}$. However, this is much more involved than for the RWM-within-Gibbs in Section A.1. Thus, it will be more convenient to construct the sets $F_d$ through the preliminary lemmas which lead to the proof of Theorems 4.1 and 4.2. The next step will involve a Taylor series expansion of $G_d V(\mathbf{x}^d)$ to show that, for large $d$, $GV(x_1)$ is a good approximation for $G_d V(\mathbf{x}^d)$. Thus, we begin by studying $\log(\frac{\pi_d(\mathbf{Y}^d)}{\pi_d(\mathbf{x}^d)})$.

LEMMA A.9. *There exists a sequence of sets $F_{d,1} \in \mathbb{R}^d$, with $\lim_{d\to\infty}\{d^{1/3}\pi_d(F^C_{d,1})\} = 0$, such that, for $\chi_i^d = 1$,*

$$\log\left\{\frac{f(Y_i^d)q(Y_i^d, x_i^d)}{f(x_i^d)q(x_i^d, Y_i^d)}\right\} = C_3(x_i^d, Z_i)d^{-1/2} + C_4(x_i^d, Z_i)d^{-2/3} + C_5(x_i^d, Z_i)d^{-5/6} + C_6(x_i^d, Z_i)d^{-1} + C_7(x_i^d, Z_i, \sigma_d),$$

*where*

$$C_3(x_i^d, Z_i) = l^3\{-\tfrac{1}{4}Z_i g'(x_i^d)g''(x_i^d) - \tfrac{1}{12}Z_i^3 g'''(x_i^d)\},$$

*and where $C_4(x_i^d, Z_i)$, $C_5(x_i^d, Z_i)$ and $C_6(x_i^d, Z_i)$ are polynomials in $Z_i$ and the derivatives of $g$. Furthermore, if $\mathbb{E}_Z$ and $\mathbb{E}_X$ denote expectation with $Z \sim N(0,1)$ and $X$ having density $f(\cdot)$, respectively, then*

$$\mathbb{E}_X\mathbb{E}_Z[C_3(X,Z)] = \mathbb{E}_X\mathbb{E}_Z[C_4(X,Z)] = \mathbb{E}_X\mathbb{E}_Z[C_5(X,Z)] = 0, \tag{A.23}$$

*whereas*

$$\mathbb{E}_X\mathbb{E}_Z[C_3(X,Z)^2] = l^6 K^2 = -2\mathbb{E}_X\mathbb{E}_Z[C_6(X,Z)]. \tag{A.24}$$

*In addition,*

$$\sup_{\mathbf{x}^d \in F_d} \mathbb{E}^*\left[\left|\sum_{i=2}^d \log\left\{\frac{f(Y_i^d)q(Y_i^d, x_i^d)}{f(x_i^d)q(x_i^d, Y_i^d)}\right\} - \left\{d^{-1/2}\sum_{i=2}^d \chi_i^d C_3(x_i^d, Z_i) - \frac{cl^6K^2}{2}\right\}\right|\right] \to 0 \qquad \text{as } d \to \infty. \tag{A.25}$$

PROOF. With the exception of (A.25) and the exact form of the sets $F_{d,1}$, the lemma is proved in [7], Lemma 1.

For $j=4,5,6$ and $x\in\mathbb{R}$, set $c_j(x)=\mathbb{E}_Z[C_j(x,Z)]$ and $v_j(x)=\mathrm{var}_Z(C_j(x,Z))$. The set $F_{d,1,j}=\bigcap_{k=1}^{3}F_{d,1,j,k}$, where

$$F_{d,1,j,1}=\left\{\mathbf{x}^d;\left|\sum_{i=2}^{d}\{C_j(x_i^d)-\mathbb{E}_X[C_j(X)]\}\right|<d^{5/8}\right\},$$

$$F_{d,1,j,2}=\left\{\mathbf{x}^d;\left|\sum_{i=2}^{d}\{V_j(x_i^d)-\mathbb{E}_X[V_j(X)]\}\right|<d^{6/5}\right\},$$

$$F_{d,1,j,3}=\left\{\mathbf{x}^d;\left|\sum_{i=2}^{d}\{C_j(x_i^d)-\mathbb{E}_X[C_j(X)]\}^2\right|<d^{6/5}\right\}.$$

Then for $j=4,5,6$ and $k=1,2,3$, it is straightforward, using Markov's inequality and conditions (4.1) and (4.2), to show that

$$d^{1/3}\pi_d(F_{d,1,j,k}^C)\to 0\qquad\text{as } d\to\infty.$$

(Cf. [7], Lemma 1, where only the cases $k=1,2$ are required.)

Finally, let $\{F_{d,1,7}\subseteq\mathbb{R}^d\}$ correspond to the sets $\{F_{n,7}\}$ constructed in [7], Lemma 1, and so, $d^{1/3}\pi_d(F_{d,1}^C)\to 0$ as $d\to\infty$, where $F_{d,1}=\bigcap_{j=4}^{7}F_{d,1,j}$.

The proof of (A.25) is then essentially the same as the proof of the final expression in [7], Lemma 1, and, hence, the details are omitted. $\square$

The next step is to find a convenient approximation for $G_dV(\mathbf{x}^d)$ which effectively allows us to consider separately the first component and the remaining $(d-1)$ components.

Lemma A.10. *Let*

$$\tilde{G}_dV(\mathbf{x}^d)=c_d d^{1/3}\mathbb{E}^*[V(Y_1)-V(x_1)]\mathbb{E}^*[e^{B_d}\wedge 1],$$

*where* $B_d(=B_d(\mathbf{x}^d))=\sum_{i=2}^{d}\{(g(Y_i^d)-g(x_i^d))-\frac{1}{\sigma_d^2}\{(x_i^d-Y_i^d-\frac{\sigma_d^2}{2}g'(Y_i^d))^2-(Y_i^d-x_i^d-\frac{\sigma_d^2}{2}g'(x_i^d))^2\}$. *There exists sets* $F_{d,2}\subseteq\mathbb{R}^d$ *with* $\lim_{d\to\infty}d^{1/3}\pi_d(F_{d,2}^C)=0$ *such that, for any* $V\in C_c^\infty$,

$$\sup_{\mathbf{x}^d\in F_{d,2}}|G_dV(\mathbf{x}^d)-\tilde{G}_dV(\mathbf{x}^d)|\to 0\qquad\text{as } d\to\infty.$$

*Moreover,*

$$\sup_{\mathbf{x}^d\in F_{d,2}}\mathbb{E}^*\left[\left|\left(\frac{\pi_d(\mathbf{Y}^d)q(\mathbf{Y}^d,\mathbf{x}^d)}{\pi_d(\mathbf{x}^d)q(\mathbf{x}^d,\mathbf{Y}^d)}\wedge 1\right)-(e^{B_d}\wedge 1)\right|\right]\to 0 \tag{A.26}$$

*as* $d\to\infty$.

PROOF. Since, conditional upon $\chi_1^d$, $Y_1$ and $\mathbf{Y}^{d-}$ are independent, it follows that

$$\tilde{G}_d V(\mathbf{x}^d) = c_d d^{1/3} \mathbb{E}^*[(V(Y_1) - V(x_1^d))(e^{B_d} \wedge 1)].$$

The lemma then follows by identical arguments to those used in [7], Theorem 3, with the sets $\{F_{d,2}\}$ chosen to correspond to the sets $\{S_n\}$ in [7], Theorem 3. $\square$

LEMMA A.11. *Let* $F_{d,3} = \{\mathbf{x}^d; g'(x_1^d) \leq d^{1/12}\}$ *then* $d^{1/3}\pi_d(F_{d,3}^C) \to 0$ *as* $d \to \infty$ *and for any* $V \in C_c^\infty$,

$$\sup_{\mathbf{x}^d \in F_{d,3}} \left| d^{1/3} c_d \mathbb{E}^*[V(Y_1^d) - V(x_1^d)] - c\frac{l^2}{2}\{g'(x_1)V'(x_1) + V''(x_1)\} \right| \to 0$$

*as* $d \to \infty$,

*with* $x_1 = x_1^d$.

PROOF. The proof is identical to [7], Lemma 2 and is, hence, omitted. $\square$

We now focus on the remaining $(d-1)$ components. First we introduce the following notation. Let $a(x) = -\frac{1}{4}g'(x)g''(x)$ and $b(x) = -\frac{1}{12}g'''(x)$. Therefore, we have that

$$C_3(x, z) = l^3\{a(x)z + b(x)z^3\}.$$

Set

$$Q_d^*(\mathbf{x}^d, \cdot) = \mathcal{L}\left\{ \frac{1}{\sqrt{d}} \sum_{i=2}^d \chi_i^d C_3(x_i^d, Z_i) \middle| \chi_1^d = 1 \right\}.$$

Let $\phi_d(\mathbf{x}^d, t) = \int_{\mathbb{R}} \exp(itw) Q_d^*(dw)$ and let $\phi(t) = \exp(-\frac{t^2}{2}cl^6K^2)$.

LEMMA A.12. *There exists a sequence of sets* $F_{d,4} \subseteq \mathbb{R}^d$ *such that:*

(a) $\lim_{d\to\infty}\{d^{1/3}\pi_d(F_{d,4}^C)\} = 0$,
(b) *for all* $t \in \mathbb{R}$,

$$\sup_{\mathbf{x}^d \in F_{d,4}} |\phi_d(\mathbf{x}^d; t) - \phi(t)| \to 0 \qquad \text{as } d \to \infty,$$

(c) *for all bounded continuous functions* $r$,

$$\sup_{\mathbf{x}^d \in F_{d,4}} \left| \int_{\mathbb{R}} Q_d^*(\mathbf{x}^d, dy) r(y) - \frac{1}{\sqrt{2\pi c l^3 K}} \int_{\mathbb{R}} r(y) \exp\left(-\frac{y^2}{cl^6K^2}\right) dy \right| \to 0$$

*as* $d \to \infty$,

(d)

$$\sup_{\mathbf{x}^d\in F_{d,4}}\left|\mathbb{E}^*\left[1\wedge\exp\left\{d^{-1/2}\sum_{i=2}^d\chi_i^d C_3(x_i^d,Z_i)-\frac{cl^6K^2}{2}\right\}\right]-2\Phi\left(-\frac{\sqrt{c}l^3K}{2}\right)\right|\to 0$$

*as* $d\to\infty$.

PROOF. The sets $F_{d,4}$ are constructed as in the proof of [7], Lemma 3, and so, statement (a) follows. Specifically, we let $F_{d,4}$ be the set of $\mathbf{x}^d\in\mathbb{R}^d$ such that

$$\left|\frac{1}{d}\sum_{i=2}^d h(x_i^d)-\int h(x)f(x)\,dx\right|\le d^{-1/4} \tag{A.27}$$

$$|h(x_i^d)|\le d^{3/4},\qquad 1\le i\le d, \tag{A.28}$$

for each of the functionals $h(x)=a(x)^2,b(x)^2,a(x)b(x),a(x)^4,b(x)^4,$ $a(x)^3b(x),a(x)^2b(x)^2,a(x)b(x)^3$.

Since statements (c) and (d) follow from statement (b) as outlined in [7], Lemma 3, all that is required is to prove (b).

Let $L_d=\{j;\chi_j^d=1,2\le j\le d\}$ and let

$$\theta_j^d(x_j^d;t)=\mathbb{E}\left[\exp\left(\frac{it}{\sqrt{d}}C_3(x_j^d,Z_j)\right)\right].$$

Let

$$\phi_d^{\Lambda_d}(\mathbf{x}^d;t)=\mathbb{E}^*\left[\exp\left\{\frac{it}{\sqrt{d}}\sum_{j\in\Lambda_d}C_3(x_j^d,Z_j)\right\}\middle|L_d=\Lambda_d\right].$$

Then since $\{C_3(x_j^d,Z_j)\}_{j=2}^d$ are independent random variables, it follows that

$$\phi_d^{\Lambda_d}(\mathbf{x}^d;t)=\prod_{j\in\Lambda_d}\theta_j^d(x_j^d;t).$$

Therefore,

$$\phi_d(\mathbf{x}^d;t)=\mathbb{E}^*\left[\exp\left\{\frac{it}{\sqrt{d}}\sum_{j=2}^d\chi_j^dC_3(x_j^d,Z_j)\right\}\right]=\mathbb{E}^*\left[\prod_{j\in L_d}\theta_j^d(x_j^d;t)\right]$$

and so,

$$\begin{aligned}&\sup_{\mathbf{x}^d\in F_{d,4}}\left|\phi_d(\mathbf{x}^d;t)-\mathbb{E}^*\left[\prod_{j\in L_d}\left\{1-\frac{t^2}{2d}v(x_j^d)\right\}\right]\right|\\&\qquad\le\sup_{\mathbf{x}^d\in F_{d,4}}\mathbb{E}^*\left[\left|\prod_{j\in L_d}\theta_j^d(x_j^d;t)-\prod_{j\in L_d}\left\{1-\frac{t^2}{2d}v(x_j^d)\right\}\right|\right],\end{aligned} \tag{A.29}$$

where $v(x_j^d) = \mathrm{var}_Z(C_3(x_j^d, Z)) = l^6\{a(x_j^d)^2 + 6a(x_j^d)b(x_j^d) + 15b(x_j^d)^2\}$.

The right-hand side of (A.29) converges to 0 as $d \to \infty$ by arguments similar to those used in [7], Lemma 3. Hence, the details are omitted.

Now by using a Taylor series expansion for $\exp(-\sum_{j=2}^d \frac{t^2}{2d}\chi_j^d v(x_j^d))$, it is trivial to show that

$$\begin{aligned}\sup_{\mathbf{x}^d \in F_{d,4}} \Bigg| \mathbb{E}^*\Bigg[\prod_{j\in L_d}\Big\{1 - \frac{t^2}{2d}v(x_j^d)\Big\}\Bigg] \\ - \mathbb{E}^*\Bigg[\exp\Bigg\{-\sum_{j=2}^d \frac{t^2}{2d}\chi_j^d v(x_j^d)\Bigg\}\Bigg]\Bigg| \to 0 \qquad \text{as } d \to \infty,\end{aligned} \tag{A.30}$$

since for all $\mathbf{x}^d \in F_{d,4}$, $\frac{1}{d^2}\sum_{j=2}^d v(x_j^d)^2 \to 0$ as $d \to \infty$ (cf. [7], Lemma 3).

The final step to complete the proof of statement (b) is to show that

$$\sup_{\mathbf{x}^d \in F_{d,4}} \Bigg| \mathbb{E}^*\Bigg[\exp\Bigg(-\sum_{j=2}^d \chi_j^d \frac{t^2}{2d} v(x_j^d)\Bigg)\Bigg] - \exp\Bigg(-\frac{t^2}{2}cl^6K^2\Bigg)\Bigg| \to 0 \qquad \text{as } d \to \infty.$$

This follows immediately, since using Chebyshev's inequality, we can show that, for all $\varepsilon > 0$,

$$\sup_{\mathbf{x}^d \in F_{d,4}} \mathbb{P}^*\Bigg(\Bigg|\sum_{j=2}^d \chi_j^d \frac{t^2}{2d} v(x_j^d) - \frac{t^2}{2}cl^6K^2\Bigg| > \varepsilon\Bigg) \to 0 \qquad \text{as } d \to \infty.$$

Thus, statement (b) is proved and the lemma follows. □

We are now in position to prove Theorems 4.1 and 4.2.

Proof of Theorem 4.1. The theorem follows from (A.25), (A.26) and part (d) of Lemma A.12. □

Proof of Theorem 4.2. We take $F_d = F_{d,1} \cap F_{d,2} \cap F_{d,3} \cap F_{d,4}$. Then

$$d^{1/3}\pi_d(F_d^C) \to 0 \qquad \text{as } d \to \infty,$$

and so, for fixed $T$,

$$\mathbb{P}(\Gamma_t^d \in F_d, 0 \le t \le T) \to 1 \qquad \text{as } d \to \infty.$$

Also, from Lemmas A.9–A.12, it follows that

$$\sup_{\mathbf{x}^d \in F_d} |G_d V(\mathbf{x}^d) - GV(x_1)| \to 0 \qquad \text{as } d \to \infty$$

for all $V \in C_c^\infty$, which depend only on the first coordinate. Therefore, the weak convergence follows by [4], Chapter 4, Corollary 8.7, since $C_c^\infty$ separates points and an identical argument to that of Theorem 3.1 can be used to demonstrate compact containment.

The maximizing of $h_c(l)$ is straightforward using the proof of [7], Theorem 2. □

**A.3. Proofs of Section 5.** The proof of Theorem 5.1 is very similar to the proof of Theorem 3.1 given in Section A.1.

First, for $\mathbf{x}^\infty \in \mathbb{R}^\infty$, let $\mathbf{x}^\infty = (x_1^\infty, x_2^\infty, \ldots)$, $\bar{x}_d = \frac{1}{d-2}\sum_{i=3}^d x_i^\infty$ and let $\bar{x} = \lim_{d\to\infty} \bar{x}_d$, should the limit exist. For $\mathbf{x}^\infty \in \mathbb{R}^\infty$, let $\mathbf{x}^d \in \mathbb{R}^d$ be such that $\mathbf{x}^d = (x_1^\infty, x_2^\infty, \ldots, x_d^\infty)$ [$= (x_1^d, x_2^d, \ldots, x_d^d)$, say], that is, $\mathbf{x}^d$ comprises the first $d$ components of $\mathbf{x}^\infty$. Then let $G_d$ be the (discrete-time) generator of $\mathbf{X}^d$, and let $V \in C_c^\infty$ be an arbitrary test function of $x_1, x_2$ and $\bar{x}_d$ only. Thus,

$$G_d V(\mathbf{x}^d) = d\mathbb{E}\biggl[(V(\mathbf{Y}^d) - V(\mathbf{x}^d))\biggl\{1 \wedge \frac{\pi_d(\mathbf{Y}^d)}{\pi_d(\mathbf{x}^d)}\biggr\}\biggr].$$

The generator $G$ of the three-dimensional diffusion described in Theorem 5.1, for an arbitrary test function $V$ of $x_1, x_2$ and $\bar{x}$, is given by

$$\begin{aligned} GV(\mathbf{x}^\infty) = \frac{cl^2}{2}\biggl(2\Phi\biggl(\frac{cl}{2\sqrt{1-\rho}}\biggr)\biggr) \\ \times \sum_{i=1}^2 \biggl\{-\frac{1}{1-\rho}(x_i - \bar{x})\frac{\partial}{\partial x_i}V(\mathbf{x}^\infty) + \frac{\partial^2}{\partial x_i^2}V(\mathbf{x}^\infty)\biggr\}. \end{aligned}$$

We shall define sets $\{F_d \subseteq \mathbb{R}^\infty; d \geq 1\}$ such that for $d\mathbb{P}(\mathbf{X}^d \in F_d^C) \to 0$ as $d \to \infty$. This is done in Lemma A.13 and, thus, we can restrict attention to $\mathbf{x}^\infty \in F_d$. Furthermore, Lemma A.13 ensures that, for all $\mathbf{x}^\infty \in F_d$, $\lim_{d\to\infty} \bar{x}_d$ exists. Therefore, since we can restrict attention to $\mathbf{x}^\infty \in F_d$, we aim to show that

$$\sup_{\mathbf{x}^\infty \in F_d} |G_d V(\mathbf{x}^d) - GV(\mathbf{x}^\infty)| \to 0 \qquad \text{as } d \to \infty, \tag{A.31}$$

which is proved in Theorem A.17 and then Theorem 5.1 follows trivially.

Then define sets $\{F_d \subseteq \mathbb{R}^\infty; d \geq 1\}$ such that for $d\mathbb{P}(\mathbf{X}^d \in F_d^C) \to 0$ as $d \to \infty$. This is done in Lemma A.13 and, thus, we can restrict attention to $\mathbf{x}^\infty \in F_d$.

LEMMA A.13. *For $1 \leq k \leq 5$, define the sequence of sets $\{F_{d,k} \subseteq \mathbb{R}^\infty; d \geq 1\}$ by*

$$F_{d,1} = \{\mathbf{x}^\infty; |R_d(\mathbf{x}^d) - (1-\rho)| < d^{-1/8}\},$$

$$F_{d,2} = \{\mathbf{x}^\infty; |\bar{x}_d - \bar{x}| < d^{-1/8}\},$$

$$F_{d,3} = \Bigl\{\mathbf{x}^\infty; \max_{1\leq i\leq d} |x_i^\infty| < d^{1/8}\Bigr\},$$

$$F_{d,4} = \biggl\{\mathbf{x}^\infty; \biggl|\frac{1}{d}\sum_{i=1}^d (x_i^\infty)^2\biggr| < d^{1/8}\biggr\},$$

$$F_{d,5} = \left\{ \mathbf{x}^\infty ; \left| \frac{1}{d} \sum_{i=1}^{d} \left( \frac{1}{1-\rho} x_i^\infty + \theta_d \sum_{j=1}^{d} x_j^\infty \right)^4 \right| < d^{1/8} \right\},$$

*where* $R_d(\mathbf{x}^d) = \frac{1}{d-1} \sum_{i=1}^d (x_i^\infty - \frac{1}{d} \sum_{j=1}^d x_j^\infty)^2$ *and* $\theta_d = -\frac{\rho}{1+(d-2)\rho-(d-1)\rho^2}$. *Let* $F_d = \bigcap_{k=1}^5 F_{d,k}$, *then*

$$d\mathbb{P}(\mathbf{X}^d \in F_d^C) \to 0 \qquad \text{as } d \to \infty. \tag{A.32}$$

PROOF. It is sufficient to show that, for $1 \le k \le 5$,

$$d\mathbb{P}(\mathbf{X}^d \in F_{d,k}^C) \to 0 \qquad \text{as } d \to \infty.$$

For the cases $k = 1, 3, 4$ and 5, it is straightforward but tedious using Markov's inequality to prove the result. Therefore, the details are omitted.

For the case $k = 2$, let $\bar{X}_d = \frac{1}{d-2} \sum_{i=3}^d X_i^\infty$ $(d \ge 3)$ and let $\bar{X} = \lim_{d\to\infty} \bar{X}_d$. Therefore, by construction (see Section 5), for all $d \ge 3$,

$$\begin{pmatrix} \bar{X}_d \\ \bar{X} \end{pmatrix} \sim N\left( \begin{pmatrix} 0 \\ 0 \end{pmatrix}, \begin{pmatrix} \frac{1}{d-2}(1+(d-3)\rho) & \rho \\ \rho & \rho \end{pmatrix} \right).$$

Thus,

$$\{\bar{X}_d | \bar{X} = \bar{x}\} \sim N\left( \frac{(d-2)\rho}{1+(d-3)\rho} \bar{x}, \frac{\rho(1-\rho)}{1+(d-3)\rho} \right). \tag{A.33}$$

Therefore, by Markov's inequality,

$$\mathbb{P}(\mathbf{X}^d \in F_{d,2}^C) = \mathbb{P}(|\bar{X}_d - \bar{X}| \ge d^{-1/8}) \le \sqrt{d}\mathbb{E}[|\bar{X}_d - \bar{X}|^4], \tag{A.34}$$

and the result follows trivially from (A.33) and (A.34). □

The procedure now differs slightly from that given in Section A.1. We postpone the finding of a suitable Taylor series expansion for $G_d V(\mathbf{x}^d)$ and first give Lemmas A.14, A.15 and A.16, which mirror Lemmas A.4, A.6 and A.7, respectively. The proofs of the aforementioned lemmas are similar to the proofs of the corresponding results in Section A.1 and, hence, the details are omitted.

LEMMA A.14. *For* $1 \le k \le d$, *let*

$$\lambda_d^k (= \lambda_d^k(\mathbf{x}^d)) = \frac{1}{d-1} \sum_{i \ne k} \chi_i^d \left( \frac{1}{1-\rho} x_i^d + \theta_d \sum_{j=1}^d x_j^d \right)^2.$$

*Then for any* $\varepsilon > 0$,

$$\sup_{\mathbf{x}^\infty \in F_d} \mathbb{P}_k^*\left( \left| \lambda_d^k - \frac{c}{1-\rho} \right| > \varepsilon \right) \to 0 \qquad \text{as } d \to \infty,$$

*where, for any random variable* $X$ *and any subset* $A \subseteq \mathbb{R}$, $\mathbb{P}_k^*(X \in A) = \mathbb{P}(X \in A | \chi_k^d = 1)$ *and* $\mathbb{E}_k^*[X] = \mathbb{E}[X | \chi_k^d = 1]$.

For $z \in \mathbb{R}$, let

$$h_d^k(z)(= h_d^k(z; \mathbf{x}^d)) = \left\{ \log\left\{ \frac{\pi_d(\mathbf{Y}^d)}{\pi_d(\mathbf{x}^d)} \right\} \Big| Z_k = z \right\}.$$

The role of $h_d^k(z)$ is similar to that played by $h_d(z)$ in Section A.1, with $h_d^k(0)$ equivalent to $B_d$ (cf. Lemma A.3).

Lemma A.15. *For $d \geq 1$ and $1 \leq k \leq d$, let $A_d^k(= A_d^k(\mathbf{x}^d)) = -\frac{cl^2}{2(1-\rho)} - \sigma_d \sum_{i \neq k} \chi_i^d (\frac{1}{1-\rho} x_i + \theta_d \sum_{j=1}^d x_j) Z_i$. Then,*

$$\sup_{\mathbf{x}^d \in F_d} |\mathbb{E}_k^*[1 \wedge e^{A_d^k}] - \mathbb{E}_k^*[1 \wedge e^{h_d^k(0)}]| \to 0 \qquad \text{as } d \to \infty, \tag{A.35}$$

*and*

$$\sup_{\mathbf{x}^d \in F_d} |\mathbb{E}_k^*[e^{A_d^k}; A_d^k < 0] - \mathbb{E}_k^*[e^{h_d^k(0)}; h_d^k(0) < 0]| \to 0 \qquad \text{as } d \to \infty.$$

Lemma A.16. *For $k \geq 1$,*

$$\sup_{\mathbf{x}^d \in F_d} \left| \mathbb{E}_k^*[1 \wedge e^{A_d^k}] - 2\Phi\left( -\frac{l}{2}\sqrt{\frac{c}{1-\rho}} \right) \right| \to 0 \qquad \text{as } d \to \infty \tag{A.36}$$

*and*

$$\sup_{\mathbf{x}^d \in F_d} \left| \mathbb{E}_k^*[e^{A_d^k}; A_d^k < 0] - \Phi\left( -\frac{l}{2}\sqrt{\frac{c}{1-\rho}} \right) \right| \to 0 \qquad \text{as } d \to \infty. \tag{A.37}$$

We are now in position to prove (A.31).

Theorem A.17.

$$\sup_{\mathbf{x}^\infty \in F_d} |G_d V(\mathbf{x}^d) - GV(\mathbf{x}^\infty)| \to 0 \qquad \text{as } d \to \infty,$$

Proof. Note that, for all $d \geq 3$, we have the following Taylor series expansion, for $V$:

$$\begin{aligned} &V(\mathbf{Y}^d) - V(\mathbf{x}^d) \\ &\quad = \sigma_d \left\{ \chi_1^d Z_1 \frac{\partial}{\partial x_1} V(\mathbf{x}^d) + \chi_2^d Z_2 \frac{\partial}{\partial x_2} V(\mathbf{x}^d) \right. \\ &\qquad \qquad \left. + \left( \frac{1}{d-2} \sum_{i=3}^d \chi_i^d Z_i \right) \frac{\partial}{\partial \bar{x}_d} V(\mathbf{x}^d) \right\} \\ &\qquad + \frac{1}{2}\sigma_d^2 \left\{ \chi_1^d Z_1^2 \frac{\partial^2}{\partial x_1^2} V(\mathbf{x}^d) + \chi_2^d Z_2^2 \frac{\partial^2}{\partial x_2^2} V(\mathbf{x}^d) + \chi_1^d \chi_2^d Z_1 Z_2 \frac{\partial^2}{\partial x_1 x_2} V(\mathbf{x}^d) \right. \end{aligned}$$

$$
\begin{aligned}
&+\chi_1^d Z_1\left(\frac{1}{d-2}\sum_{i=3}^{d}\chi_i^d Z_i\right)\frac{\partial^2}{\partial x_1\,\partial\bar{x}_d}V(\mathbf{x}^d)\\
&+\chi_2^d Z_2\left(\frac{1}{d-2}\sum_{i=3}^{d}\chi_i^d Z_i\right)\frac{\partial^2}{\partial x_2\,\partial\bar{x}_d}V(\mathbf{x}^d)\\
&\qquad+\left(\frac{1}{d-2}\sum_{i=3}^{d}\chi_i^d Z_i\right)^2\frac{\partial^2}{\partial\bar{x}_d^2}V(\mathbf{x}^d)\Bigg\}\\
&+\frac{1}{6}\sigma_d^3 F(\mathbf{x}^d,\chi^d,\mathbf{Z})\\
&=\sum_{i=1}^{9}D_i(\mathbf{x}^d,\chi^d,\mathbf{Z})+\frac{1}{6}\sigma_d^3 F(\mathbf{x}^d,\chi^d,\mathbf{Z}),\qquad\text{say,}
\end{aligned}
$$

where $F(\mathbf{x}^d,\chi^d,\mathbf{Z})$ is a function of $\chi^d=(\chi_1^d,\chi_2^d,\dots,\chi_d^d)$, $\mathbf{Z}$ and the third derivatives of $V(\mathbf{x}^d)$. Since $V\in C_c^\infty$, it follows that, for all $\mathbf{x}^d\in\mathbb{R}^d$, $\mathbb{E}[|F(\mathbf{x}^d,\chi^d,\mathbf{Z})|]<\infty$, and so,

$$
\sup_{\mathbf{x}^d\in F_d}|G_dV(\mathbf{x}^d)-\hat{G}_dV(\mathbf{x}^d)|\to 0\qquad\text{as } d\to\infty,
$$

where

$$
\begin{aligned}
\hat{G}_dV(\mathbf{x}^d)&=\sum_{i=1}^{9}\mathbb{E}\left[D_i(\mathbf{x}^d,\chi^d,\mathbf{Z})\left\{1\wedge\frac{\pi(\mathbf{Y}^d)}{\pi(\mathbf{x}^d)}\right\}\right]\\
&=\sum_{i=1}^{9}\hat{G}_d^iV(\mathbf{x}^d),\qquad\text{say.}
\end{aligned}
$$

Now for all $\mathbf{x}^d\in F_d$, we have that

$$
\begin{aligned}
\hat{G}_d^1V(\mathbf{x}^d)&=d\sigma_d\mathbb{E}\left[\chi_1^d Z_1\frac{\partial}{\partial x_1}V(\mathbf{x}^d)\{1\wedge\exp(h_d^1(Z_1))\}\right]\\
&=dc_d\sigma_d\mathbb{E}_1^*\Bigg[Z_1\frac{\partial}{\partial x_1}V(\mathbf{x}^d)\Bigg\{\{1\wedge\exp(h_d^1(0))\}\\
&\qquad-\sigma_d Z_1\left(\frac{1}{1-\rho}x_1+\theta_d\sum_{j=1}^{d}x_j^d\right)\\
&\qquad\times\{\exp(h_d^1(0));h_d^1(0)<0\}\Bigg\}\Bigg]\\
&\quad+O(d\sigma_d^3 d^{1/4}).
\end{aligned}
$$

Therefore, since $Z_1$ and $h_d^1(0)$ are independent, it follows that

$$\sup_{\mathbf{x}^d\in F_d}\left|\hat{G}_d^1 V(\mathbf{x}^d) - dc_d\sigma_d^2\left\{-\left(\frac{1}{1-\rho}x_1+\theta_d\sum_{j=1}^d x_j^d\right)\right\}\right.$$

$$\left.\times\frac{\partial}{\partial x_1}V(\mathbf{x}^d)\mathbb{E}_1^*[\exp(h_d^1(0)); h_d^1(0)<0]\right|\to 0 \tag{A.38}$$

$$\text{as } d\to\infty.$$

Now for all $\mathbf{x}^d\in F_d$,

$$\begin{aligned}\theta_d\sum_{i=1}^d x_i &= \theta_d(x_1+x_2) - \frac{\rho(d-2)}{1+(d-2)\rho-(d-1)\rho^2}\left\{\frac{1}{d-2}\sum_{i=3}^d x_i\right\}\\ &\to -\frac{1}{1-\rho}\bar{x}\qquad\text{as } d\to\infty.\end{aligned}\tag{A.39}$$

Therefore, it follows from (A.38), (A.39) and Lemma A.16 that

$$\sup_{\mathbf{x}^d\in F_d}\left|\hat{G}_d^1 V(\mathbf{x}^d) - \frac{cl^2}{2(1-\rho)}\left\{2\Phi\left(-\frac{l}{2}\sqrt{\frac{c}{1-\rho}}\right)\right\}\right.$$

$$\left.\times\left\{-(x_1-\bar{x})\frac{\partial}{\partial x_1}V(\mathbf{x}^\infty)\right\}\right|\to 0\qquad\text{as } d\to\infty.$$

Similarly,

$$\sup_{\mathbf{x}^d\in F_d}\left|\hat{G}_d^2 V(\mathbf{x}^d) - \frac{cl^2}{2(1-\rho)}\left\{2\Phi\left(-\frac{l}{2}\sqrt{\frac{c}{1-\rho}}\right)\right\}\right.$$

$$\left.\times\left\{-(x_2-\bar{x})\frac{\partial}{\partial x_2}V(\mathbf{x}^\infty)\right\}\right|\to 0\qquad\text{as } d\to\infty.$$

Next, for all $\mathbf{x}^d\in F_d$, we have that

$$\begin{aligned}\hat{G}_d^3 V(\mathbf{x}^d) &= d\sigma_d\mathbb{E}\left[\frac{\partial}{\partial\bar{x}_d}V(\mathbf{x})\left(\frac{1}{d-2}\sum_{i=3}^d\chi_i^d Z_i\right)\left\{1\wedge\frac{\pi(\mathbf{Y}^d)}{\pi(\mathbf{x}^d)}\right\}\right]\\ &= dc_d\sigma_d\frac{\partial}{\partial\bar{x}_d}V(\mathbf{x}^d)\\ &\quad\times\left\{\frac{1}{d-2}\sum_{i=3}^d\mathbb{E}_i^*\left[Z_i\left\{\{1\wedge\exp(h_d^i(0))\}\right.\right.\right.\\ &\qquad\qquad\left.\left.\left. -\sigma_d Z_i\left(\frac{1}{1-\rho}x_i+\theta_d\sum_{j=1}^d x_j\right)\right.\right.\right.\end{aligned}$$

$$\times \{\exp(h_d^i(0)); h_d^i(0) < 0\}\Bigg\}\Bigg]\Bigg\}$$

$$+ O(d\sigma_d^3 d^{1/4}).$$

Therefore, since $Z_i$ and $h_d^i(0)$ are independent, it follows that

$$\sup_{\mathbf{x}^d \in F_d} |\hat{G}_d^3 V(\mathbf{x}^d) - \check{G}_d^3 V(\mathbf{x}^d)| \to 0 \qquad \text{as } d \to \infty,$$

where

$$\check{G}_d^3 V(\mathbf{x}^d) = dc_d\sigma_d^2 \frac{\partial}{\partial \bar{x}_d} V(\mathbf{x}^d)\Bigg\{-\frac{1}{d-2}\sum_{i=3}^d \Bigg(\frac{1}{1-\rho}x_i + \theta_d \sum_{j=1}^d x_j\Bigg)$$

$$\times \mathbb{E}_i^*[\exp(h_d^i(0)); h_d^i(0) < 0]\Bigg\}.$$

Let

$$\tilde{G}_d^3 V(\mathbf{x}^d) = \frac{cl^2}{1-\rho}\frac{\partial}{\partial \bar{x}_d} V(\mathbf{x}^d)\Bigg\{-\frac{1}{d-2}\sum_{i=3}^d \mathbb{E}_i^*[\exp(h_d^i(0)); h_d^i(0) < 0](x_i - \bar{x})\Bigg\}.$$

Then since, for all $\mathbf{x}^d \in F_d$, $d\theta_d \to -\frac{1}{1-\rho}$ and $\frac{1}{d}\sum_{j=1}^d x_j^d \to \bar{x}$, as $d \to \infty$, we have that

$$\sup_{\mathbf{x}^d \in F_d} |\check{G}_d^3 V(\mathbf{x}^d) - \tilde{G}_d^3 V(\mathbf{x}^d)| \to 0 \qquad \text{as } d \to \infty.$$

By Lemma A.16, for all $\mathbf{x}^d \in F_d$ and $i \geq 3$,

$$\mathbb{E}_k^*[\exp(h_d^i(0)); h_d^i(0) < 0] \to \Phi\Bigg(-\frac{l}{2}\sqrt{\frac{c}{1-\rho}}\Bigg) \qquad \text{as } d \to \infty.$$

Therefore, since $\bar{x}_d \to \bar{x}$ as $d \to \infty$, we have that

$$\sup_{\mathbf{x}^d \in F_d} |\tilde{G}_d^3 V(\mathbf{x}^d)| \to 0 \qquad \text{as } d \to \infty.$$

Hence,

$$\sup_{\mathbf{x}^d \in F_d} |\hat{G}_d^3 V(\mathbf{x}^d)| \to 0 \qquad \text{as } d \to \infty.$$

Now for all $\mathbf{x}^d \in F_d$, we have, by independence, that

$$\begin{aligned}\hat{G}_d^4 V(\mathbf{x}) &= \frac{1}{2} d\sigma_d^2 c_d \mathbb{E}_1^*\Bigg[Z_1^2 \frac{\partial^2}{\partial x_1^2} V(\mathbf{x}^d)\{1 \wedge \exp(h_d^1(Z_1))\}\Bigg] \\ &= \frac{1}{2} d\sigma_d^2 c_d \mathbb{E}_1^*\Bigg[Z_1^2 \frac{\partial^2}{\partial x_1^2} V(\mathbf{x}^d)\{1 \wedge \exp(h_d^1(0))\}\Bigg] + O(d\sigma_d^3 d^{1/8}) \\ &= \frac{1}{2} d\sigma_d^2 c_d \frac{\partial^2}{\partial x_1^2} V(\mathbf{x}^d)\mathbb{E}_1^*[1 \wedge \exp(h_d^1(0))] + O(d\sigma_d^3 d^{1/8}),\end{aligned}$$

since, for all $\mathbf{x}^d \in F_d$, $|\frac{1}{1-\rho}x_1 + \theta_d \sum_{i=1}^d x_i^d| \le d^{1/8}$. Therefore, by Lemma A.16, we have that

$$\sup_{\mathbf{x}^d \in F_d} \left| \hat{G}_d^4 V(\mathbf{x}^d) - \frac{cl}{2}\left\{2\Phi\left(\frac{l}{2}\sqrt{\frac{c}{1-\rho}}\right)\right\}\frac{\partial^2}{\partial x_1^2}V(\mathbf{x}^\infty)\right| \to 0 \qquad \text{as } d \to \infty.$$

Similarly,

$$\sup_{\mathbf{x}^d \in F_d} \left| \hat{G}_d^5 V(\mathbf{x}^d) - \frac{cl}{2}\left\{2\Phi\left(\frac{l}{2}\sqrt{\frac{c}{1-\rho}}\right)\right\}\frac{\partial^2}{\partial x_2^2}V(\mathbf{x}^\infty)\right| \to 0 \qquad \text{as } d \to \infty.$$

There exists $\alpha_1$ lying between 0 and $Z_1$, such that

$$\begin{aligned}
&|\hat{G}_d^6 V(\mathbf{x}^d)| \\
&\quad = \left|\frac{1}{2}\sigma_d^2 \mathbb{E}\left[\chi_1^d \chi_2^d Z_1 Z_2 \frac{\partial^2}{\partial x_1\,\partial x_2}V(\mathbf{x}^d)\left\{1 \wedge \frac{\pi(\mathbf{Y}^d)}{\pi(\mathbf{x}^d)}\right\}\right]\right| \\
&\quad = \left|\frac{1}{2}dc_d\sigma_d^2\right. \\
&\qquad \times \mathbb{E}_1^*\left[\chi_2^d Z_2 Z_1 \frac{\partial^2}{\partial x_1\,\partial x_2}V(\mathbf{x}^d)\left\{\{1 \wedge \exp(h_d^1(0))\}\right.\right. \\
&\qquad\qquad - \sigma_d Z_1\left(\frac{1}{1-\rho}x_1 + \theta_d \sum_{j=1}^d x_j^d\right) \\
&\qquad\qquad \times \left.\left.\left.\{\exp(h_d^1(\alpha_1)); h_d^1(\alpha_1) < 0\}\right\}\right]\right|.
\end{aligned}$$

Note that $Z_1$ is independent of $\chi_2^d$, $Z_2$ and $h_d^1(0)$. Therefore, since $\mathbb{E}_1^*[Z_1] = 0$, we have that

$$|\hat{G}_d^6 V(\mathbf{x}^d)| \le \frac{1}{2}d\sigma_d^3 c_d \frac{\partial^2}{\partial x_1\,\partial x_2}V(\mathbf{x}^d)\left|\frac{1}{1-\rho}x_1 + \theta_d \sum_{j=1}^d x_j^d\right|. \tag{A.40}$$

Since $\frac{\partial^2}{\partial x_1\,\partial x_2}V(\mathbf{x}^d)$ is bounded and for all $\mathbf{x}^d \in F_d$, $|\frac{1}{1-\rho}x_1 + \theta_d \sum_{j=1}^d x_j| < d^{1/8}$, it follows that the right-hand side of (A.40) converges to 0 as $d \to \infty$. Hence, $|\hat{G}_d^6 V(\mathbf{x}^d)| \to 0$ as $d \to \infty$.

Similarly, for $i = 7, 8, 9$, it can be shown that $|\hat{G}_d^i V(\mathbf{x}^d)| \to 0$ as $d \to \infty$ and the lemma follows. □

Proof of Theorem 5.1. The proof is similar to that of Theorem 3.1. From Lemma A.13 and Theorem A.17, we have that $d\mathbb{P}(\mathbf{X}^d \notin F_d) \to 0$ as $d \to \infty$ and

$$\sup_{\mathbf{x}^\infty \in F_d} |G_d V(\mathbf{x}^d) - GV(\mathbf{x}^\infty)| \to 0 \qquad \text{as } d \to \infty,$$

respectively. Therefore, the weak convergence follows by [4], Chapter 4, Corollary 8.7, since $C_c^\infty$ separates points and a similar argument to that of Theorem 3.1 can be used to demonstrate compact containment. □

The proof of Theorem 5.3 is similar to the proof of Theorem 5.1 and, hence, the details are omitted. The key point is to show that Lemma A.13 still holds with (A.32) replaced by

$$d^2 \mathbb{P}(\mathbf{X}^d \in F_d^C) \to 0 \qquad \text{as } d \to \infty.$$

This is again straightforward, but tedious, using Markov's inequality.

School of Mathematics
University of Manchester
Sackville Street
Manchester, M60 1QD
United Kingdom
E-mail: P.Neal-2@manchester.ac.uk

Department of Mathematics and Statistics
Fylde College
Lancaster University
Lancaster, LA1 4YF
United Kingdom
E-mail: g.o.roberts@lancaster.ac.uk